# ASYMPTOTIC APPROXIMATION OF NONPARAMETRIC REGRESSION EXPERIMENTS WITH UNKNOWN VARIANCES[1]

By Andrew V. Carter

*University of California, Santa Barbara*

Asymptotic equivalence results for nonparametric regression experiments have always assumed that the variances of the observations are known. In practice, however the variance of each observation is generally considered to be an unknown nuisance parameter. We establish an asymptotic approximation to the nonparametric regression experiment when the value of the variance is an additional parameter to be estimated or tested. This asymptotically equivalent experiment has two components: the first contains all the information about the variance and the second has all the information about the mean. The result can be extended to regression problems where the variance varies slowly from observation to observation.

**1. Introduction.** We will show that a nonparametric regression experiment where the variance is unknown (and possibly changing) is asymptotically equivalent to a continuous Gaussian process. This equivalence is demonstrated by the explicit construction of the continuous Gaussian process from the nonparametric regression observations and vice versa.

In particular, a simple version of the nonparametric regression problem observes $n$ independent normals,

$$(1.1) \qquad Y_i = f(i/n) + \sigma \xi_i, \qquad i = 1, \ldots, n,$$

where $f$ is an unknown smooth function that we want to estimate (or test), the $\xi_i$ are independent standard normals and $\sigma^2$ is the variance of the noise.

Brown and Low [2] showed that this nonparametric regression problem is asymptotically equivalent to trying to estimate $f$ in the white-noise experiment that observes the continuous process

$$(1.2) \qquad Y(t) = \int_0^t f(x)\,dx + \frac{\sigma}{\sqrt{n}} W(t), \qquad 0 \le t \le 1,$$

Received September 2005; revised November 2006.

[1]Supported by NSF Grant DMS-05-04233.

where $W(t)$ is a standard Brownian motion (SBM). Brown and Low [2] assumed that the variance structure was known, and their construction of the white-noise process $Y(t)$ depends crucially on the value of $\sigma$. In practice, however, we typically do not know the value of $\sigma$, and it is usually considered a secondary "nuisance" parameter. Its estimation is only necessary to the extent that it calibrates the estimation of $f$ (as in setting a threshold level or bandwidth). Our approach is to include the variance as a second parameter so that the experiment is now concerned with decision procedures concerning the pair $(f, \sigma)$. While it is too strong to assume that $\sigma$ is known, we may be erring in the other direction by promoting its importance to the same level as the mean function. However, Theorem 1 shows that there is no significant penalty to pay in treating the variance as part of the parameter space because the equivalence holds for essentially the same spaces as in [2].

Our motivation comes from wavelet thresholding techniques (e.g., [5, 7]) that estimate the variance using the high frequency wavelet coefficients and estimate the mean mainly from the low frequency coefficients using the estimate of the variance to determine which terms to include in the model. A similar approach is used by Rice [21] to choose the bandwidths of kernel estimators. Our asymptotic approximation contains two components: a $\chi^2$-distributed random variable with information about the variance, and a continuous Gaussian process with information about the mean.

This sort of approximation is also available when the variance is a function over the unit interval, $Y_i = f(i/n) + \sigma(i/n)\xi_i$. In this case, the strategy is to separate the $Y_i$'s into groups such that within each group the variance function is nearly constant and then to proceed as in the constant variance case. Not surprisingly, if the variance is also to be nonparametrically estimated, the equivalence result is only true under somewhat stricter conditions on the means.

### 1.1. *Asymptotic equivalence.*

The proposed approximation is in the sense of Le Cam's deficiency distance between statistical experiments [15]. This type of approximation provides a correspondence between estimation procedures in each experiment such that any good estimator in the asymptotic approximation corresponds to a good estimator in the nonparametric regression estimator and vice versa.

In this formulation, we have a pair of statistical experiments $\mathcal{P}$ and $\mathcal{Q}$ that consist of sets of distribution functions $\{\mathbb{P}_f \mid f \in \mathcal{F}\}$ on $(\mathcal{X}, \mathcal{A})$ and $\{\mathbb{Q}_f \mid f \in \mathcal{F}\}$ on $(\mathcal{Y}, \mathcal{B})$. Both are indexed by the same parameter set $\mathcal{F}$, and the two experiments are equivalent if they provide the same information about $f \in \mathcal{F}$. Le Cam proposed a pseudometric for statistical experiments $\Delta(\mathcal{P}, \mathcal{Q}) = \max[\delta(\mathcal{P}, \mathcal{Q}), \delta(\mathcal{Q}, \mathcal{P})]$ where $\delta(\mathcal{P}, \mathcal{Q}) = \inf_K \sup_f \|K(\mathbb{P}_f) - \mathbb{Q}_f\|_{\mathrm{TV}}$, using the total variation distance and "transitions" $K$ that map distributions on the

sample space of $\mathcal{P}$ to distributions on the sample space of $\mathcal{Q}$. For our purposes, however, Le Cam's general notion of transitions (see [16], page 18) is not necessary. Instead, we will bound $\delta(\mathcal{P}, \mathcal{Q})$ by proposing a randomized transformation of the observed data. Thus, $K$ can be represented by the conditional distribution on $(\mathcal{Y}, \mathcal{B})$ given an observation from $\mathbb{P}_f$, and $K(\mathbb{P}_f)$ is the marginal distribution on $\mathcal{Y}$.

Therefore, the first step in bounding this $\Delta$-distance is to propose a candidate transformation from $\mathcal{X}$ to $\mathcal{Y}$. Then the bound is established by bounding the distance between the distributions of the transformed observations and the observations from the approximating experiment. Explicit transformations between the experiments are useful because they generate a correspondence between the estimators in experiments. For instance, if the distribution of $T(X)$ is close to that of $Y$, then the estimator $\hat{f}(T(X))$ has nearly the same risk as $\hat{f}(Y)$. The transformation $T$ may be randomized in the sense that it may depend on some external random variables, but it may not depend on the parameters. In particular, the transformation in [2] depends on the variance $\sigma^2$ which is now a part of our parameter space. Therefore, we must formulate a different transformation that does not depend on $\sigma^2$.

Two sequences of experiments $\mathcal{Q}_n$ and $\mathcal{P}_n$ are asymptotically equivalent if $\Delta(\mathcal{P}_n, \mathcal{Q}_n) \to 0$. Asymptotic equivalence implies that the risk under a bounded loss function achieved by any estimator in $\mathcal{P}_n$ can be achieved asymptotically by associated estimators in $\mathcal{Q}_n$ and vice versa [15].

### 1.2. Main results.

1.2.1. *Constant variance.* In order to accommodate the added aspect of an unknown variance, the parameter space will be expanded to include both the smooth functions $f$ and the variance $\sigma^2$. The specification of the exact set of parameters is described in Section 2. For Theorem 1, the parameter space $\mathcal{F}_\sigma$ includes all functions in a Hölder space with $\alpha > 1/2$ as in [2].

THEOREM 1. *Suppose that the experiment $\mathcal{P}_n$ observes $Y_i$ as in* (1.1). *The distributions are indexed by* $(f, \sigma) \in \mathcal{F}_\sigma \times \mathbb{R}^+$ *as in Definition* 1.

*Further, suppose that the experiment $\mathcal{Q}_n$ has distributions indexed by the same pairs* $(f, \sigma)$ *with* $V \sim \Gamma(\frac{n}{2}, \frac{2\sigma^2}{n})$ *and*

$$(1.3) \qquad Y(t) \mid V = \int_0^t f(x)\,dx + V^{1/2} n^{-1/2} W(t), \qquad 0 \le t \le 1,$$

*where $W(t)$ is a SBM.*

*Then the experiments $\mathcal{P}_n$ and $\mathcal{Q}_n$ are asymptotically equivalent,* $\Delta(\mathcal{P}_n, \mathcal{Q}_n) \to 0$.

The proof for simplified versions of these experiments is in Section 3, and the rest of the argument is in Section 5.

*Remarks on Theorem* 1. Note that the experiment $\mathcal{Q}_n$ is equivalent to observing $Y(t)$ alone because the random variable $V$ can be computed almost surely from $Y(t)$ via its quadratic variation. This is why the variance of $Y(t)$ is random as opposed to just $\sigma^2/n$ as in [2]. $\mathcal{Q}_n$ is like Le Cam's locally asymptotically mixed normal experiment ([17], page 121).

One implication of this approximation is that asymptotically $V$ is sufficient for estimating $\sigma^2$. Further, because the distributions of $Y(t)$ conditional on different values of $V$ are mutually singular, the conditionality principle implies that inference for $f$ should be performed conditional on $V$. For instance, minimax results for the white-noise problem similar to those in [20] or [6] bound $\inf_{\hat{f}} \sup_{f \in \mathcal{F}_\sigma} \mathbb{E}_f L(\hat{f}, f) = \mathcal{R}(\mathcal{F}_\sigma, \sigma^2)$ for observations as in (1.2). This implies a bound on the risk in estimating $f$ in (1.3) using the expected value of the conditional risk, $\mathbb{E}\mathcal{R}(\mathcal{F}, V)$. Theorem 1 implies that the same asymptotic minimax result also applies to $\mathcal{P}_n$ (for bounded loss functions).

The advantage of approximating the regression observations in (1.1) by the continuous process is that it makes certain calculations easier. For instance, in the experiment $\mathcal{Q}_n$ a linear estimate of the mean $f$ at a point $t$ which is of the form $Y(K_h)$, where $K_h$ is a kernel function with bandwidth $h$, has a normal distribution with mean $\int K_h(x) f(x)\,dx$ and variance $\frac{V}{n}\int K_h^2(x)\,dx$ conditional on $V$. The bandwidth $h$ should be chosen as a function of $V$ to minimize the error, and even knowing the exact value of $\sigma$ would not provide a better bandwidth. This is very much like the approach to choosing a bandwidth described in [21].

1.2.2. *Varying variance.* It may be more interesting to consider nonparametric regression experiments where the variance changes over the interval. Here the parameter space is the product of two function spaces, $\breve{\mathcal{F}}_\sigma$ and $\Sigma$. The varying variance means that more smoothness is required in $\breve{\mathcal{F}}_\sigma$, in particular $\breve{\mathcal{F}}_\sigma$ includes Hölder spaces only for $\alpha > 3/4$ and $\Sigma$ includes functions $\sigma(t)$ such that $\log \sigma(t)$ is Hölder for $\alpha > 1$. The exact definition of the parameter space is in Section 2.

THEOREM 2. *The experiment* $\breve{\mathcal{P}}_n$ *observes* $Y_i = f(i/n) + \sigma(i/n)\xi_i$ *for* $i = 1, \ldots, n$ *with* $\xi_i$ *as in Theorem* 1 *and* $(f, \sigma^2) \in \breve{\mathcal{F}}_\sigma(\alpha, \gamma_k) \times \Sigma(\alpha_1)$, *the parameter space in Definition* 2.

*The experiment* $\breve{\mathcal{Q}}_{n,m}$ *observes two Gaussian processes:*

$$V(t) = \int_0^t \log \sigma^2(x)\,dx + \sqrt{2}n^{-1/2}W_2(t),$$

*and then conditional on* $V(t)$,

$$dY(t) = f(t)\,dt + Z_\ell n^{-1/2}\,dW_1(t), \qquad (\ell-1)/m < t \le \ell/m,$$

as $\ell = 1, \ldots, m$, where $Z_\ell = \exp[\frac{m}{2}(V(\ell/m) - V([\ell-1]/m))]$. The processes $W_1(t)$ and $W_2(t)$ are independent SBMs.

For any $\alpha > 3/4$ and $\alpha_1 > \max(1, \frac{\alpha}{2\alpha-1})$, there is a sequence $m_n$ (where $n^{1/3} < m_n < n^{1/2}$ depending on $\alpha$ and $\alpha_1$) such that these experiments are asymptotically equivalent, $\Delta(\breve{\mathcal{P}}_n, \breve{\mathcal{Q}}_{n,m_n}) \to 0$.

The proof of a simplified version of this result is in Section 4 and the rest of the proof is in Section 6 and depends on asymptotic results in Sections 7 and 8 that bound the distance between the simplified experiments and $\breve{\mathcal{P}}$ and $\breve{\mathcal{Q}}$, respectively.

Theorem 2 requires a bit more smoothness on $f$, but we can trade off a bit of smoothness in the set of mean functions for less smoothness in the variance space. In particular, if $f$ is a Hölder function for some $3/4 < \alpha < 1$, then $\log \sigma^2(x)$ needs to be a Hölder function for $\alpha_1 > 2\alpha/(2\alpha - 1)$. This always implies that $\alpha_1 > 1$ and the variance function has one bounded derivative.

Sections 9 and 10 give bounds on the K–L divergence between gamma and normal distributions that will be used in the proofs. Further technical lemmas are established in Sections 11, 12 and 13.

### 1.3. *Related work.*

The equivalence results of Brown and Low [2] were the first to apply Le Cam's deficiency to a nonparametric regression experiment, and their introduction provides a number of further references motivating this approach. Brown, Cai, Low and Zhang [1] extends their results to the case where the design points are randomly chosen uniformly over the interval. Rohde [22] uses a Fourier series decomposition to get better results in the approximation of Brown and Low [2]. Carter [3] extends the fixed design result to the unit square. All of these results assume that the errors are normal with a known variance. Brown and Low [2] discuss how to adjust their results to cases where the variance changes over the interval, but it is essential to their methodology that the experimenter know the variance at each observation.

Grama and Nussbaum [10] discusses nonparametric regression problems with nonnormal errors. In particular, one of the cases they treat is estimating the variance of normal observations. Zhou [24] treats the variance case in particular and improves the bounds to apply to Besov spaces. The $\tilde{\mathcal{Q}}$ experiment in Theorem 2 reduces to the continuous Gaussian experiment from [10] if the mean function is assumed known. Therefore, Theorem 2 synthesizes the results of both [2] and [10] (under stronger smoothness conditions).

The most interesting applications of our work may be in heteroscedastic nonparametric regression, of which there is a considerable literature, (e.g., [4, 9, 11, 13, 23]). The variance estimator of Müller and Stadtmuller [19] seems closest to these results in that the mean squared error is estimated on a fine grid, thus producing approximately the $V(t)$ or $Z_\ell$ of experiment $\check{Q}_{n,m}$, and then these observations are smoothed to produce an estimator of $\sigma^2(t)$. Hall, Kay and Titterington [11] improve on this technique by finding the best linear functions of the $Y_i$ which can be squared and averaged to get a estimate of $\sigma^2$. Other types of estimators from Fan and Yao [9] and Ruppert, Wand, Holst and Hössjer [23] are based on residuals from a preliminary fit of the mean. Testing for heteroscedasticity is addressed in Dette and Munk [4] based on differences between successive observations, while Eubank and Thomas [8] use residuals.

**2. The parameter spaces.** The most convenient way of describing these parameter spaces is using wavelet bases and their associated Besov sequence norms. Assuming that the mean functions $f$ are in $L^2([0,1))$, let $\phi_{k,j}$ for $j = 1, \ldots, 2^k$ and $\psi_{i,j}$ for $i \geq k$ and $j = 1, \ldots, 2^i$ be the scale functions and wavelets, respectively, for an orthonormal wavelet basis on $[0, 1]$. For most of our arguments it is necessary that these are the Haar basis: $\phi_{k,j}(t) = 2^{k/2}\phi_0(2^k t - j + 1)$ where $\phi_0(t) = \mathbf{1}\{0 \leq t \leq 1\}$, and $\psi_{i,j}$ is defined analogously with $\psi_0(t) = \mathbf{1}\{0 \leq t \leq \frac{1}{2}\} - \mathbf{1}\{\frac{1}{2} \leq t \leq 1\}$.

The coefficients $\vartheta_{k,j}$ and $\theta_{i,j}$ are such that

$$(2.1) \qquad f(x) = \sum_{j=1}^{2^k} \vartheta_{k,j}\phi_{k,j}(x) + \sum_{i=k}^{\infty}\sum_{j=1}^{2^i}\theta_{i,j}\psi_{i,j}(x).$$

These coefficients can be found via $\theta_{i,j} = \langle f, \psi_{i,j} \rangle = \int f \psi_{i,j}\, dx$. The parameter space is the set of smooth functions that can be described succinctly by the basis functions $\phi_{k,j}$ and $\psi_{i,j}$.

There are two Besov sequence norms we will use on the series of coefficients. The $b(\alpha, 2, 2)$ norm is

$$\|\boldsymbol{\theta}\|_{b(\alpha,2,2)} = \left(\sum_j \vartheta_j^2 + \sum_{i \geq k} 2^{2\alpha i}\sum_j \theta_{i,j}^2\right)^{1/2}$$

and the $b(\alpha, \infty, 1)$ norm is

$$\|\boldsymbol{\theta}\|_{b(\alpha,\infty,1)} = \sup_j 2^{k(\alpha+1/2)}|\vartheta_j| + \sum_{i \geq k} 2^{i(\alpha+1/2)}\sup_j |\theta_{i,j}|.$$

These norms are equivalent to Besov norms (see, e.g., [12], Chapter 9).

The parameter spaces we will use are compact in these norms. This implies that there is a uniform bound on the partial sums in each norm. Specifically, for a sequence of positive $\gamma_k$ that goes to 0 as $k$ gets large, let $\Theta(\alpha, 2, 2, \gamma_k)$ be the set of all sequences $\theta_{i,j}$ such that

$$(2.2) \qquad \sum_{i \geq k} 2^{2\alpha i}\sum_j \theta_{i,j}^2 \leq \gamma_k^2.$$

Analogously,

$$\sup_{\theta \in \Theta(\alpha, \infty, 1, \gamma_k)} \left[ \sum_{i \geq k} 2^{i(\alpha+1/2)} \sup_j |\theta_{i,j}| \right] \leq \gamma_k.$$

The results in Sections 3 and 4 require that the mean functions are in $\Theta(\alpha, 2, \gamma_k)$ while the approximations in Sections 5 and 6 require that the means are in $\Theta(\alpha, \infty, 1, \gamma_k)$. Hölder$(M, \alpha^*)$ functions with $\alpha < \alpha^* < 1$ are in both spaces with $\gamma_k = M 2^{k(\alpha - \alpha^*)}$.

DEFINITION 1. Using the Haar basis functions, the parameter space for Theorem 1 is

$$\mathcal{F}_\sigma(\gamma_k) \times \mathbb{R}^+ = \{(f, \sigma^2) : f \in \Theta(1/2, 2, 2, \gamma_k \sigma) \cap \Theta(1/2, \infty, 1, \gamma_k \sigma), \sigma^2 > 0\}.$$

When the constant variance is replaced by a function $\sigma^2(t)$, then greater smoothness in the mean function is necessary. The log of the variance function is assumed to be a Hölder function with $\alpha_1 > 1$. Let $\tau(t) = \log \sigma^2(t)$. Then the parameter set H$(M, \alpha_1)$ includes all such $\tau$ where

$$(2.3) \quad \sup_{t \in [0,1]} |\tau'(t)| \leq M \quad \text{and} \quad \sup_{t,s \in [0,1]} |\tau'(t) - \tau'(s)| \leq M|s - t|^{\alpha_1 - 1}.$$

Furthermore, let $\log \bar{\sigma}^2 = \int_0^1 \log \sigma^2(t) \, dt$.

DEFINITION 2. The parameter space for Theorem 2 is

$$\breve{\mathcal{F}}_\sigma(\alpha, \gamma_k) \times \Sigma(\alpha_1) = \{(f, \sigma^2) : f(t) \in \Theta(\alpha, 2, 2, \gamma_k \bar{\sigma}) \cap \Theta(\alpha, \infty, 1, \gamma_k \bar{\sigma}),$$
$$\text{and } \log(\sigma^2(t)) \in H(M, \alpha_1)\},$$

where again the basis is assumed to be the Haar basis.

REMARK. There are two tricks used here to avoid the condition $\sigma^2(t) > \varepsilon > 0$. First the smoothness of the functions is measured relative to the variance in that the tail of the Besov norm has to decrease proportionally with $\sigma$. Also, the smoothness condition on the variance is on the logarithm of $\sigma^2(t)$ as opposed to $\sigma^2(t)$ itself, thus $\sup_{0 < t < 1} |\log \frac{\sigma^2(t)}{\bar{\sigma}^2}| \leq M$ by the mean value theorem. We avoided a lower bound on the variance so that the experiments will still be invariant under rescalings.

**3. Sequence space result.** Instead of working directly $\mathcal{P}_n$, we will first consider an experiment based on the orthonormal basis functions. The experiment $\bar{\mathcal{P}}_k$ observes $n = 2^k$ independent normals

$$X_{0,j} \sim \mathcal{N}\left(\vartheta_{k_0,j}, \frac{\sigma^2}{n}\right) \qquad \text{for } j = 1, \ldots, 2^{k_0},$$

$$X_{i,j} \sim \mathcal{N}\left(\theta_{i,j}, \frac{\sigma^2}{n}\right) \qquad \text{for } k_0 \leq i < k,$$

where the $\vartheta_{k_0,j}$ and $\theta_{i,j}$ are the wavelet coefficients of the mean function $f$ as defined above.

The experiment $\bar{\mathcal{P}}_\infty$ observes the entire sequence of normal random variables. This sequence experiment is equivalent to the experiment that observes the Gaussian process from (1.2) because the random coefficients can be generated from this process via $X_{0,j} = \int \phi_{k_0,j}(t)\,dY(t)$ and $X_{i,j} = \int \psi_{i,j}(t)\,dY(t)$, and the process can be constructed from the coefficients via (2.1) using $\theta_{i,j} = X_{i,j}$. Unfortunately, this experiment $\bar{\mathcal{P}}_\infty$ is completely informative with respect to estimating $\sigma^2$ and therefore cannot be asymptotically equivalent to $\mathcal{P}_n$ except in trivial cases.

Instead, $\bar{\mathcal{P}}_k$ is approximated by $\mathcal{Q}$ which replaces $\sigma^2$ in $\bar{\mathcal{P}}_\infty$ with a chi-squared observation. The experiment $\mathcal{Q}$ observes the random variables

$$V \sim \Gamma\left(\frac{n}{2}, \frac{2\sigma^2}{n}\right), \qquad Y_{0,j}|V \sim \mathcal{N}\left(\vartheta_{k_0,j}, \frac{V}{n}\right), \qquad Y_{i,j}|V \sim \mathcal{N}\left(\theta_{i,j}, \frac{V}{n}\right).$$

LEMMA 1. *For the parameter set* $\bar{\mathcal{F}}_\sigma(\gamma_k) \times \mathbb{R}^+ = \{(\boldsymbol{\theta}, \sigma^2) : \boldsymbol{\theta} \in \Theta(1/2, 2, 2, \gamma_k\sigma), \sigma^2 \in \mathbb{R}^+\}$ *and* $n = 2^k$, $m = 2^{k_0}$ *and* $m = n\gamma_{k_0}$,

$$\Delta(\bar{\mathcal{P}}_k, \mathcal{Q}) \leq 2\gamma_{k_0}^{1/2}. \tag{3.1}$$

*If* $\gamma_k = M2^{-\varepsilon k}$ *for some small* $\varepsilon$ *between 0 and 1/2, then*

$$\Delta(\bar{\mathcal{P}}_k, \mathcal{Q}) \leq 2M^{1/2}n^{-\varepsilon/(2(1+\varepsilon))}. \tag{3.2}$$

Clearly, (3.2) follows directly from (3.1) using $m = M^{1/(1+\varepsilon)}n^{1/(1+\varepsilon)}$.

This lemma and its proof imply that there is a sense in which the $\chi^2$ random variable $V$ and the scaling function coefficients $Y_{0,j}$ are asymptotically sufficient statistics for these experiments. Sections 3.5 and 3.7 argue that the information about $f$ in the $Y_{i,j}$ and $X_{i,j}$ is negligible and the information about the variances is summarized in $V$.

The bound on the deficiency $\delta(\bar{\mathcal{P}}_k, \mathcal{Q})$ is described in Sections 3.2–3.6 and the bound on $\delta(\mathcal{Q}, \bar{\mathcal{P}}_k)$ is in Section 3.7.

### 3.1. *Using Kullback–Leibler divergence.*

The total variation distance measures the distance between distributions in the deficiency distance, but total variation is inconvenient for the analysis especially in the case of product measures. It is easier to establish the bounds using Kullback–Leibler divergence because the divergence between the joint distributions of $X$ and $Y$ is equal to the divergence between the marginal distributions of $X$ plus the expected value of the divergence between the conditional distributions of $Y$ given $X$. In particular, for product measures $\mathbf{D}(\prod_i P_i, \prod_i Q_i) = \sum_i \mathbf{D}(P_i, Q_i)$. We also have the bound $\|\mathbb{P} - \mathbb{Q}\|_{\mathrm{TV}} \leq \sqrt{\mathbf{D}(\mathbb{P}, \mathbb{Q})}$ which allows us to apply K–L bounds to the deficiency distance.

In a convenient abuse of notation, we will use $\mathbf{D}(X, Y)$ to refer to the divergence between the *distributions* of $X$ and $Y$. This allows us to write

$$\mathbf{D}((X_1, Y_1,), (X_2, Y_2)) = \mathbf{D}(X_1, X_2) + \mathbb{E}[\mathbf{D}(Y_1, Y_2 \mid X)].$$

Furthermore, $\mathbf{D}(X, Y) \leq \mathbb{E}[\mathbf{D}(X, Y \mid V)]$, which means that the divergence between the marginal distributions of $X$ and $Y$ is less than the expected value of the divergence between the conditional distribution of $X$ given $V$ and the conditional distribution of $Y$ given $V$.

3.2. *Hierarchical structure of the experiments.* The strategy is to consider the $X_{0,j}$ and $Y_{0,j}$ observations as those containing the information about the mean, and the rest as having information about the variance. The $\mathcal{Q}$ experiment consists of observations of $V$, $Y_{0,j}$ and $Y_{i,j}$ for $i \geq k_0$ as above.

We can construct a parallel structure from the $\bar{\mathcal{P}}_k$ observations,

$$\hat{V} = \frac{n}{n-m} \sum_{k_0 < i \leq k} \sum_j X_{i,j}^2,$$

$$\hat{Y}_{0,j} = X_{0,j} \mid \hat{V} \sim \mathcal{N}\left(\vartheta_{k_0,j}, \frac{\sigma^2}{n}\right), \qquad \hat{Y}_{i,j} \mid \hat{V} \sim \mathcal{N}\left(0, \frac{\hat{V}}{n}\right),$$

where the $\hat{Y}_{i,j}$ are generated independently conditional on the estimated variance.

In both experiments, the conditional distribution given $V$ of the wavelet coefficients, $Y_{i,j}$, is independent of the distribution of the scaling function coefficients $Y_{0,j}$. Thus, there is a decomposition of the bound on the divergence into three terms,

$$\begin{aligned}
(3.3) \qquad \mathbf{D}(\mathbb{Q}, K(\mathbb{P})) = {} & \mathbf{D}(V, \hat{V}) + \sum_j \mathbb{E}[\mathbf{D}(Y_{0,j}, X_{0,j} \mid V)] \\
& + \sum_{i \geq k_0} \sum_j \mathbb{E}[\mathbf{D}(Y_{i,j}, \hat{Y}_{i,j} \mid V)].
\end{aligned}$$

We will bound the contribution from each of these terms in Sections 3.3, 3.4 and 3.5, respectively.

3.3. *The variances.* The construction first generates an estimate of the variance. This estimate, $\hat{V}$, is $\sigma^2/(n-m)$ times a noncentral $\chi^2$ random variable with $n - m$ degrees of freedom and noncentrality parameter $\mu = n \sum_{k_0 \leq i < k} \sum_j \theta_{i,j}^2 / \sigma^2$.

The distributions of $\hat{V}$ and $V$ are approximately gamma with $\alpha = (n - m)/2$ and $n/2$, respectively. The distribution of $\hat{V}$ is approximate because the normals have nonzero means. Ignoring that inconvenience for a moment, Lemma 6 bounds the divergence by $\mathbf{D}(V, \hat{V}) \approx m^2 n^{-2}$.

The effect of the noncentrality of $\hat{V}$ on the bound can be handled using a mixture distribution characterization of the noncentral $\chi^2$. A noncentral $\chi_n^2$ with noncentrality parameter $\mu$ can be generated from a Poisson mixture of $\chi^2$ distributions with $2\Lambda + n$ degrees of freedom with mixing parameter $\Lambda \sim \text{Poisson}(\mu)$. For the K–L divergence between mixture distributions, $\mathbf{D}(V, \hat{V}) \le \mathbb{E}_\mu \mathbf{D}(V, \hat{V} \mid \Lambda)$ where $V$ is independent of $\Lambda$. From (9.6),

$$\mathbb{E}_\mu \mathbf{D}(V, \hat{V} \mid \Lambda) \le \frac{m^2}{2n^2} + \mathbb{E}_\mu\left(\frac{\Lambda^2}{2(n-m)} + \frac{\Lambda}{n-1}\right) \le \frac{m^2}{2n^2} + \frac{\mu^2 + 3\mu}{2(n-m)}.$$

The size of $\mu$ will be shown below in (3.7) to be $n\gamma_{k_0}^2/m$, and therefore only the first term in this bound will concern us,

$$(3.4) \qquad \mathbf{D}(V, \hat{V}) \le \frac{m^2}{2n^2} + \text{smaller order terms.}$$

3.4. *The top level.* The broadest coefficients $X_{0,j}$ are equated directly to the $Y_{0,j}$, and the distributions are normals with means $\hat{\vartheta}_{k_0,j}$. The only difference between the two sets of coefficients is the variance: for the $Y_{0,j}$ it is $V/n$, and for the $X_{0,j}$ it is $\sigma^2/n$.

Therefore, from (10.1),

$$\sum_j \mathbb{E}[\mathbf{D}(Y_{0,j}, X_{0,j} \mid V)] = \frac{m}{2}\mathbb{E}\left[\frac{V}{\sigma^2} - 1 - \log\left(\frac{V}{\sigma^2}\right)\right]$$

$$(3.5) \hspace{4cm} \le \frac{m}{2}\log\left(\frac{n}{n-2}\right)$$

$$\le \frac{m}{n-2}$$

using Jensen's inequality on $\mathbb{E}\log 1/V$. This is actually larger than the error in (3.4), and they both imply that $m = o(n)$.

3.5. *The bottom levels.* The third and final term in (3.3) compares the wavelet coefficients. The coefficients $\hat{Y}_{i,j}$ are uninformatively generated by random zero-mean normals with variance $\hat{V}/n$.

The difference between one of these normals and a normal generated by the experiment $\mathcal{Q}$ (conditional on $V$) is in the difference of the means, $\mathbf{D}(Y_{i,j}, \hat{Y}_{i,j} \mid V) = n\theta_{i,j}^2 V^{-1}/2$. Thus the total error is

$$(3.6) \quad \sum_{i \ge k_0}\sum_j \mathbb{E}[\mathbf{D}(Y_{i,j}, \hat{Y}_{i,j} \mid V)] = \frac{n^2}{2(n-2)}\sum_{i,j}\frac{\theta_{i,j}^2}{\sigma^2} \le n\sum_{i>k_0}\sum_j \frac{\theta_{i,j}^2}{\sigma^2}.$$

Using $m = 2^{k_0}$,

$$(3.7) \qquad n\sum_{i>k_0}\sum_j \frac{\theta_{i,j}^2}{\sigma^2} \le \frac{n}{m}\sum_{i>k_0}2^i\sum_j \frac{\theta_{i,j}^2}{\sigma^2} \le \frac{n}{m}\gamma_{k_0}^2.$$

3.6. *Choosing* $m$. Choosing the dimension of the scaling functions balances the errors in the approximation of the scaling and the wavelet coefficients. The trade-off is between the bound in (3.5) ($m/n$) and the bound in (3.7). Minimizing the bound is then possible by setting the two terms equal to each other,

$$\frac{m}{n} = \frac{n\gamma_{k_0}^2}{m} \quad \Longrightarrow \quad m = n\gamma_{k_0} \quad \Longrightarrow \quad \frac{m^2}{n^2} + \frac{m}{n} + \frac{n}{m}\gamma_{k_0}^2 \le 3\gamma_{k_0}.$$

Therefore, plugging (3.4), (3.5) and (3.7) into (3.3) for $m = n\gamma_{k_0}$ gives $\delta(\bar{\mathcal{P}}_k, \mathcal{Q}) \le 2\gamma_{k_0}^{1/2}$.

3.7. *The transition in the other direction.* The other half of the deficiency distance bound requires a map from the $(V, Y_{i,j})$ observations from $\mathcal{Q}$ to the $X_{i,j}$ observations from $\bar{\mathcal{P}}_k$. Once again, the top level observations remain unchanged, $\hat{X}_{0,j} = Y_{0,j} \sim \mathcal{N}(\vartheta_{k_0,j}, V/n)$. The $Y_{i,j}$ for $i \ge k_0$ are not used in the transformation; instead the $\hat{X}_{i,j}$ are functions of the variance $V$. Because $V$ is a sufficient statistic for estimating $\sigma^2$ from $n$ independent normals with mean 0, there is a probability distribution conditional on $V$ that is not a function of $\sigma$ of $n$ independent $\mathcal{N}(0, \sigma^2/n)$ random variables from which we can use the first $n - m$ as $\hat{X}_{i,j}$.

The $\hat{X}_{i,j}$ and $\hat{X}_{0,j}$ are not independent because they both depend on $V$, but we can bound the K–L divergence via

$$\mathbf{D}(\hat{\mathbf{X}}, \mathbf{X}) = \sum_{k < i < k_0} \sum_j \mathbf{D}(\hat{X}_{i,j}, X_{i,j}) + \mathbb{E} \sum_j \mathbf{D}(\hat{X}_{0,j}, X_{0,j} \mid \{X_{i,j}\}_{k < i < k_0, j})$$

$$\le \sum_{k < i < k_0} \sum_j \mathbf{D}(\hat{X}_{i,j}, X_{i,j}) + \mathbb{E} \sum_j \mathbf{D}(Y_{0,j}, X_{0,j} \mid V).$$

The contribution to the error from the second term is just as in (3.5), and the first term is less than the error in (3.7). Therefore, $\delta(\mathcal{Q}, \bar{\mathcal{P}}_k) \le 2\gamma_{k_0}^{1/2}$ and Lemma 1 is established.

## 4. A variance function over the interval.
An interesting extension of the result in Lemma 1 is to consider what happens when the variance changes over the interval. Our simplified version of this experiment assumes that we can group the basis functions into $m_1 = 2^{k_1}$ groups for $k_1 < k_0$,

$$\mathcal{I}_\ell = \left\{ (i,j) : 2^i\left(\frac{\ell - 1}{m_1}\right) < j \le 2^i\left(\frac{\ell}{m_1}\right) \right\},$$

and then each group of coefficients will have a different variance, $\sigma_\ell^2$. These groups are chosen so that each Haar basis function $\psi_{i,j}$ with $(i,j) \in \mathcal{I}_\ell$ has support in $(\ell-1)/m_1 < t \le \ell/m_1$. The variance of the group will be determined by the variance function $\sigma^2(t)$ via

$$\log \sigma_\ell^2 = m_1 \int_{(\ell-1)/m_1}^{\ell/m_1} \log \sigma^2(t) \, dt.$$

LEMMA 2. *The parameter space* $\tilde{\mathcal{F}}_{\bar{\sigma}} \times \Sigma$ *contains the mean functions* $f(t) \in \Theta(\alpha, 2, 2, \gamma_k \bar{\sigma})$ *for some* $\alpha > 3/4$, *and the variance functions* $\sigma^2(t)$ *are such that* $\sup_\Sigma \max_\ell (\log \sigma_\ell^2 - \log \bar{\sigma}^2) \le M$.

*The experiment* $\tilde{\mathcal{P}}_k$ *observes* $n$ *independent normals,*

$$X_{0,j} \sim \mathcal{N}\left(\vartheta_{k_0,j}, \frac{\sigma_\ell^2}{n}\right) \qquad \textit{for } (k_0, j) \in \mathcal{I}_\ell,$$

$$X_{i,j} \sim \mathcal{N}\left(\theta_{i,j}, \frac{\sigma_\ell^2}{n}\right) \qquad \textit{for } (i, j) \in \mathcal{I}_\ell, \ k_0 \le i < k.$$

*The* $\tilde{\mathcal{Q}}$ *experiment observes* $m_1$ *independent* $V_\ell \sim \Gamma(\frac{n}{2m_1}, \frac{2\sigma_\ell^2 m_1}{n})$ *and*

$$Y_{0,j}|\mathbf{V} \sim \mathcal{N}\left(\vartheta_{k_0,j}, \frac{V_\ell}{n}\right) \qquad \textit{for } (k_0, j) \in \mathcal{I}_\ell,$$

$$Y_{i,j}|\mathbf{V} \sim \mathcal{N}\left(\theta_{i,j}, \frac{V_\ell}{n}\right) \qquad \textit{for } (i, j) \in \mathcal{I}_\ell, i \ge k_0,$$

*where the normals are all conditionally independent.*

*Then for* $m_0 = 2^{k_0}$,

$$\Delta(\tilde{\mathcal{P}}_k, \tilde{\mathcal{Q}}) \le 2m_1^{1/2} m_0^{1/2} n^{-1/2} + e^{M/2} m_0^{-\alpha} n^{1/2} \gamma_{k_0}.$$

This bound is of order $\gamma_{k_0}$ when $m_0 = n^{1/(2\alpha)}$ and $m_1 = n^{1-1/(2\alpha)} \gamma_{k_0}^2$.

The basic idea is that the argument for Lemma 1 can be repeated on each of the $m_1$ independent pieces of this experiment. A Haar basis is used here because the basis functions have disjoint support which keeps things tidy as the variance changes over the interval.

Comparing this to Lemma 1, the approximate sufficient statistics are the $m_0$ random variables $Y_{0,j}$ and the $m_1$ variances $V_\ell$ where $m_0 < m$ from Lemma 1 because the $f$ are smoother and $m_1 < m_0$ because the variance functions are smoother still.

4.1. *Proof of Lemma* 2. The transformation of the $X_{i,j}$ follows as in Section 3 on each of the $m_1$ pieces. First, there are estimates of the variances based on the observations for $i \ge k_0$,

$$\hat{V}_\ell = \frac{nm_1}{n - m_1} \sum_{(i,j) \in \mathcal{I}_\ell, i \ge k_0} X_{i,j}^2.$$

Then new Gaussian observations $\hat{Y}_{i,j}$ for $i \geq k_0$ are generated by independent normals with variance $\hat{V}_\ell$.

The error in the approximation is bounded in the same three stages: first the error in the estimation of the variance, then the difference between the distributions when $i = 0$, and finally the distance between the distributions of the observations for $i \geq k_0$,

$$\begin{aligned}
\mathbf{D}(\tilde{\mathbb{Q}}, \tilde{\mathbb{P}}K) &= \sum_{\ell=1}^{m} \mathbf{D}(V_\ell, \hat{V}_\ell) + \sum_{\ell=1}^{m_1} \sum_{\{j\,:\,(k_0,j)\in\mathcal{I}_\ell\}} \mathbb{E}[\mathbf{D}(Y_{0,j}, X_{0,j} \mid V_\ell)] \\
&\quad + \sum_{\ell=1}^{m_1} \sum_{i \geq k_0} \sum_{\{j\,:\,(i,j)\in\mathcal{I}_\ell\}} \mathbb{E}[\mathbf{D}(Y_{i,j}, \hat{Y}_{i,j} \mid V_\ell)].
\end{aligned}$$

(4.1)

Each of the estimates $\hat{V}_\ell$ is $m_1\sigma_\ell^2/(n - m_1)$ times a noncentral $\chi^2$ with $(n - m_1)/m_1$ degrees of freedom and noncentrality parameter

$$\mu_\ell = \sum_{k_0 \leq i < k} \sum_{\{j\,:\,(i,j)\in\mathcal{I}_\ell\}} \frac{n\theta_{i,j}^2}{\sigma_\ell^2},$$

which is small for large $k_0$. By (9.6), the divergence between the distributions of $\hat{V}_\ell$ and $V_\ell$ is

$$\mathbf{D}(\hat{V}_\ell, V_\ell) \leq \frac{m_1^2}{2(n - m_1)^2} + \frac{m_1(\mu_\ell^2 + 2\mu_\ell)}{n}.$$

Using independence, the divergence between the distributions of the entire vectors is just the sum,

$$\sum_{\ell=1}^{m} \mathbf{D}(\hat{V}_\ell, V_\ell) \leq \frac{m_1^3}{n^2} + \frac{3m_1}{n}\left(\sum_{\ell=1}^{m} \mu_\ell\right).$$

(4.2)

These terms will turn out to be negligible relative to the errors in the other two terms.

The second term in (4.1) is the divergence between the conditional distributions of $X_{0,j}$ and the $Y_{0,j}$. By (10.1) and Jensen's inequality,

$$\mathbb{E}[\mathbf{D}(Y_{0,j}, X_{0,j} \mid V_\ell)] = \frac{1}{2}\left[\frac{\mathbb{E}(V_\ell)}{\sigma_\ell^2} - 1 - \mathbb{E}\log\left(\frac{V_\ell}{\sigma_\ell^2}\right)\right] \leq \frac{m_1}{n - 2m_1}.$$

Thus the bound on the sum is

$$\sum_{\ell=1}^{m_1} \sum_{\{j\,:\,(k_0,j)\in\mathcal{I}_\ell\}} \mathbb{E}[\mathbf{D}(Y_{0,j}, X_{0,j} \mid V_\ell)] \leq \frac{m_1 m_0}{n}.$$

(4.3)

Finally, the third term in (4.1) is the divergence between the conditional distributions of the $Y_{i,j}$'s and the $\hat{Y}_{i,j}$'s. By (10.1)

$$\sum_{\ell=1}^{m_1} \sum_{i \geq k_0} \sum_{\{j\,:\,(i,j) \in \mathcal{I}_\ell\}} \mathbb{E}[\mathbf{D}(Y_{i,j}, \hat{Y}_{i,j} \mid V_\ell)] = \sum_{\ell=1}^{m_1} \sum_{i \geq k_0} \sum_{\{j\,:\,(i,j) \in \mathcal{I}_\ell\}} \mathbb{E}\frac{n\theta_{i,j}^2}{2V_\ell}.$$

In this case, $\mathbb{E}V_\ell^{-1} = \sigma_\ell^{-2}(1 - \frac{2m_1}{n})^{-1} \leq e^M \bar{\sigma}^{-2}$ by the smoothness condition in Lemma 2. Thus

$$\sum_{\ell=1}^{m_1} \sum_{i \geq k_0} \sum_{\{j\,:\,(i,j) \in \mathcal{I}_\ell\}} \mathbb{E}[\mathbf{D}(Y_{i,j}, \hat{Y}_{i,j} \mid V_\ell)] \leq e^M n \sum_{i \geq k_0} \sum_j \frac{\theta_{i,j}^2}{\bar{\sigma}^2}.$$

Here we will use the smoothness properties of the function space to get the bound

$$(4.4) \qquad n \sum_{i \geq k_0} \sum_j \frac{\theta_{i,j}^2}{\bar{\sigma}^2} \leq n m_0^{-2\alpha} \sum_{i \geq k_0} 2^{2\alpha i} \sum_j \frac{\theta_{i,j}^2}{\bar{\sigma}^2} \leq n m_0^{-2\alpha} \gamma_{k_0}^2.$$

Therefore, plugging in (4.2), (4.3) and (4.4) into (4.1) yields

$$\mathbf{D}(\tilde{\mathbb{Q}}, \tilde{\mathbb{P}}K) \leq \frac{2m_1 m_0}{n} + e^M \frac{n}{m_0^{2\alpha}} \gamma_{k_0}^2,$$

and the deficiency is

$$\delta(\tilde{\mathcal{P}}, \tilde{\mathcal{Q}}) \leq 2m_1^{1/2} m_0^{1/2} n^{-1/2} + e^{M/2} n^{1/2} m_0^{-\alpha} \gamma_{k_0}.$$

4.1.1. *The transformation in the other direction.* The proof of Lemma 2 is completed by bounding the deficiency in the other direction. Following very much what we did in Section 3.7, each $V_\ell$ can be decomposed into $n/m_1$ independent normals with mean 0 and variance $\sigma_\ell^2$ to create a set of observations (conditional on $V_\ell$)

$$\hat{X}_{0,j} \sim \mathcal{N}\left(\vartheta_{k_0,j}, \frac{V_\ell}{n}\right), \qquad \hat{X}_{i,j} \sim \mathcal{N}\left(0, \frac{\sigma_\ell^2}{n}\right)$$

for $(i,j) \in \mathcal{I}_\ell$ and $k_0 \leq i < k$.

The divergence between these distributions and the distributions in $\tilde{\mathcal{P}}_k$ is less than

$$\sum_{\ell=1}^{m_1} \sum_{k < i \leq k_0} \sum_{\{j\,:\,(i,j) \in \mathcal{I}_\ell\}} \mathbf{D}(\hat{X}_{i,j}, X_{i,j}) + \sum_{\ell=1}^{m_1} \sum_{\{j\,:\,(k_0,j) \in \mathcal{I}_\ell\}} \mathbb{E}\mathbf{D}(\hat{X}_{0,j}, X_{0,j} \mid V_\ell),$$

where the first term is bounded as in (4.4), and the second term as in (4.3).

Therefore,

$$\delta(\tilde{\mathcal{P}}_k, \tilde{\mathcal{Q}}) \leq m_1^{1/2} m_0^{1/2} n^{-1/2} + e^{M/2} n^{1/2} m_0^{-\alpha} \gamma_{k_0}$$

and the proof of Lemma 2 is finished.

**5. The difference between the sequence and sampled experiments.** The sequence results in Lemma 1 and 2 assume that the observations are $n$ wavelet coefficients. This is unrealistic, and we would prefer that the observations were of the form (1.1). The standard technique is to use $Y_\ell / \sqrt{n}$ as approximations to the scaling coefficients at the lowest level in the wavelet expansion $\vartheta_{k,\ell}$.

The cascade algorithm of [18] constructs higher frequency scaling function coefficients from the scale function and wavelet coefficients. This construction can be used to generate $Y_\ell^* \sim \mathcal{N}(n^{1/2}\theta_{k,\ell}, \sigma^2)$ from the $Y_{i,j}$, and the construction can be inverted to construct the wavelet coefficients from these scaling function coefficients at level $k$.

The error in this approximation is in the difference between the means,

$$\sum_\ell \mathbf{D}(Y_\ell^*, Y_\ell) = \sum_{\ell=1}^n \frac{(n^{1/2}\vartheta_{k,\ell} - f(\ell/n))^2}{2\sigma^2}.$$

To be concrete, take the orthonormal basis to be the Haar basis. For $f$ a continuous function,

$$(5.1) \qquad f(\ell/n) = n^{1/2}\vartheta_{k,\ell} + \sum_{i > \log_2 n} \theta_{i,j^*} 2^{(i-1)/2},$$

where the $j^*$'s are the indices of the wavelets such that $|\phi_{i,j^*}(\ell/n)| > 0$.

The K–L divergence bound becomes

$$
\begin{aligned}
(5.2) \qquad \sum_\ell \mathbf{D}(Y_\ell^*, Y_\ell) &= \frac{1}{2\sigma^2} \sum_{\ell=1}^n \left( \sum_{i > \log_2 n} 2^{(i-1)/2} \theta_{i,j^*} \right)^2 \\
&\leq \frac{n}{2\sigma^2} \left( \sum_{i > \log_2 n} 2^{(i-1)/2} \sup_j |\theta_{i,j}| \right)^2 \\
&\leq \frac{1}{4} \left( \sum_{i > \log_2 n} 2^i \sup_j \left| \frac{\theta_{i,j}}{\sigma} \right| \right)^2 \leq \frac{1}{4}\gamma_k^2
\end{aligned}
$$

for mean functions in $\Theta(1/2, \infty, 1, \gamma_k\sigma)$. This function space is different than required in Lemma 1 but still includes any Hölder space for $\alpha > 1/2$ as in [2].

Theorem 1 is proven by first appealing to the triangle inequality for $\Delta$, $\Delta(\mathcal{P}, \mathcal{Q}) \leq \Delta(\mathcal{P}, \bar{\mathcal{P}}) + \Delta(\bar{\mathcal{P}}, \mathcal{Q})$. Then, for the parameter space that is the intersection of the two spaces required for Lemma 1 and (5.2), $\Delta(\mathcal{P}, \bar{\mathcal{P}}) \leq \gamma_k$ from (5.2) and $\Delta(\bar{\mathcal{P}}, \mathcal{Q}) \leq 2M^{1/2}\gamma_{k_0}$ by Lemma 1. Therefore, $\Delta(\mathcal{P}, \mathcal{Q}) \leq 3M^{1/2}\gamma_{k_0} \to 0$, and Theorem 1 is established.

**6. Proving Theorem 2.** There are three asymptotic results that can be combined to establish Theorem 2. The first is Lemma 2, which showed

$$\Delta(\tilde{\mathcal{P}}_k, \tilde{\mathcal{Q}}) \leq e^{M/2} n^{1/2} m_0^{-\alpha} \gamma_{k_0} + 2m_1^{1/2} m_0^{1/2} n^{-1/2}.$$

We need two more approximations. First the $\mathcal{P}$ experiment needs to be approximated by the nonparametric regression experiment.

LEMMA 3.  *We have*

$$\Delta(\tilde{\mathcal{P}}_k, \check{\mathcal{P}}) \leq 2e^{M/2}(n^{1/2} m_0^{-\alpha} \gamma_{k_0} + n^{1/2} m_1^{-\alpha_1}$$
$$+ m_0^{1/2} m_1^{-1} + n^{1/2} m_0^{-1} + n^{1/2} m_1^{-3/2}).$$

This is proven in Section 7.

Finally, we need to approximate $\tilde{\mathcal{Q}}$ by a continuous Gaussian process.

LEMMA 4.

$$\Delta(\check{Q}, \tilde{Q}) \leq m_1 n^{-1/2} + 2M n^{1/2} m_1^{-\alpha_1} + M n^{1/2} m_1^{-3/2}.$$

This result is shown in Section 8.

These three results can be combined using the triangle inequality for the $\Delta$ distance to prove that

$$\Delta(\check{\mathcal{P}}, \check{\mathcal{Q}}) \leq C n^{1/2} m_0^{-\alpha} \gamma_{k_0} + C n^{1/2} m_1^{-\alpha_1} + C m_1^{1/2} m_0^{1/2} n^{-1/2} + m_0^{1/2} m_1^{-1} + \cdots.$$

Let $\zeta_0 = \log_n m_0$ and $\zeta_1 = \log_n m_1$ so that if

$$\zeta_0 \geq \frac{1}{2\alpha}, \qquad \zeta_1 > \frac{1}{2\alpha_1}, \qquad \zeta_0 + \zeta_1 < 1 \quad \text{and} \quad \zeta_0 < 2\zeta_1,$$

then $\Delta(\check{\mathcal{P}}, \check{\mathcal{Q}}) \to 0$. The conditions can only be fulfilled if $\alpha > 3/4$ and $\alpha_1 > \alpha/(2\alpha - 1)$. Of course the argument assumes all along that $\alpha \leq 1$ and $\alpha_1 > 1$.

For $\frac{3}{4} < \alpha < 1$ and $1 < \alpha_1 < \frac{3}{2}$, we could take $m_0 = n^{1/(2\alpha)}$ and $m_1 = n^{1/(2\alpha_1)} n^{\varepsilon}$, where $\varepsilon = \frac{1}{4}(\frac{2\alpha - 1}{\alpha} - \frac{1}{\alpha_1})$ so that

$$n^{1/2} m_0^{-\alpha} \gamma_{k_0} + n^{1/2} m_1^{-\alpha_1} + m_1^{1/2} m_0^{1/2} n^{-1/2} + m_0^{1/2} m_1^{-1}$$
$$= \gamma_{k_0} + n^{-\varepsilon \alpha_1} + n^{-\varepsilon/2} + n^{1/(4\alpha) - 1/(2\alpha_1)} n^{-\varepsilon},$$

which goes to 0 as $n \to \infty$ because $\varepsilon > 0$ and $4\alpha > 2\alpha_1$. The other terms in the bound on $\Delta(\check{\mathcal{P}}, \check{\mathcal{Q}})$ are also negligible,

$$n^{1/2} m_0^{-1} + n^{1/2} m_1^{-3/2} = n^{(\alpha-1)/(2\alpha)} + n^{(2\alpha_1 - 3)/(4\alpha_1)} n^{-3\varepsilon/2},$$

which goes to 0 because $\alpha < 1$ and $\alpha_1 < \frac{3}{2}$. Finally, if $n^{1/2} m_0^{-1} \to 0$ and $m_0 m_1 n^{-1} \to 0$ then clearly $m_1 n^{-1/2} \to 0$.

Therefore, Lemmas 2, 3 and 4 together are sufficient to prove Theorem 2 where the sequence $m_n = n^{1/(2\alpha_1) + \varepsilon}$.

**7. Sequence space result for changing variances.** Lemma 3 compares the experiment $\breve{\mathcal{P}}_k$ which observes $n$ wavelet coefficients to the experiment $\breve{\mathcal{P}}$ which observes $n$ normals with means $f(i/n)$ and variances $\sigma^2(i/n)$. The approximation will be established in a series of steps by establishing intermediary experiments that are equivalent to both experiments.

7.1. *Negligible wavelet means.* The first approximation of $\breve{\mathcal{P}}$ is by $\mathcal{P}_1^*$, which observes the $X_{0,j}$ the same as in $\breve{\mathcal{P}}$ but observes $X_{i,j}^*$ that have expectation zero. The divergence between the joint distributions is

$$\mathbf{D}((\{X_{0,j}\}, \{X_{i,j}\}), (\{X_{0,j}\}, \{X_{i,j}^*\})) = n \sum_{i \geq k_0} \sum_j \frac{\theta_{i,j}^2}{\sigma_\ell^2},$$

which is the same as the bound in (4.4), and therefore it goes to 0 for $n$ large and

$$(7.1) \qquad \Delta(\breve{\mathcal{P}}_k, \mathcal{P}_{1,n}^*) \leq e^{M/2} n^{1/2} m_0^{-\alpha} \gamma_{k_0}.$$

A sequence of zero-mean normals with an unknown variance has the sum of the squared observations as a sufficient statistic. In particular, $\mathcal{P}_1^*$ is equivalent to observing

$$X_{0,j} \sim \mathcal{N}\left(\vartheta_{k_0,j}, \frac{\sigma_\ell^2}{n}\right), \qquad \hat{V}_\ell \sim \Gamma\left(\frac{n-m_0}{2m_1}, \frac{2m_1\sigma_\ell^2}{n-m_0}\right)$$

for $\ell = 1, \ldots, m_1$, where

$$\hat{V}_\ell = \frac{nm_1}{n-m_0} \sum_{(i,j) \in \mathcal{I}_\ell} X_{i,j}^{*2}.$$

7.2. *Distributing the variances.* Next, we consider the experiment $\mathcal{P}_2^*$ that has $m_0$ variances instead of $m_1$. It has observations

$$X_{0,j}^* \sim \mathcal{N}\left(\vartheta_{k_0,j}, \frac{\bar{\sigma}_j^2}{n}\right) \quad \text{and} \quad V_j^* \sim \Gamma\left(\frac{n-m_0}{2m_0}, \frac{2m_0\bar{\sigma}_j^2}{n-m_0}\right),$$

all independent for $j = 1, \ldots, m_0$. The new variances are

$$\log \bar{\sigma}_j^2 = [\zeta_{\ell,j} \log \sigma_\ell^2 + \zeta_{\ell+1,j} \log \sigma_{\ell+1}^2]$$

for $(2\ell-1)/2m_1 < j/m_0 \leq (2\ell+1)/2m_1$ where $\zeta_{\ell,j} + \zeta_{\ell+1,j} = 1$ are weight functions defined below.

This experiment is generated by smoothing out the variance information in the $\hat{V}_\ell$ of $\mathcal{P}_1^*$. The transformation of the observations from $\mathcal{P}_1^*$ leaves the observations $X_{0,j}$ alone and uses the $\hat{V}_\ell$ to generate the $\chi^2$ random variables. The trick is to redistribute this information in a smooth way to produce the $V_j^*$.

7.2.1. *The transformation of the $\hat{V}_\ell$.* We decompose each $\hat{V}_\ell$ into $2m_1/m_0$ gamma observations (with a small correction for $\ell = 1$ and $\ell = m_1$). The technique depends on the fact that a random variable $X \sim \Gamma(\alpha, \beta)$ times a beta random variable $\xi \sim \mathcal{B}(\delta\alpha, [1-\delta]\alpha)$ is $X\xi \sim \Gamma(\delta\alpha, \beta)$ and is independent of $X(1-\xi) \sim \Gamma([1-\delta]\alpha, \beta)$. This can be extended to a multivariate beta distribution. For $j$ between $(2\ell - 3)m_0/(2m_1)$ and $(2\ell + 1)m_0/(2m_1)$ and parameters $\delta_{\ell,j}$ such that $\sum_j \delta_{\ell,j} = 1$, the density of $\xi_{\ell,j}$ on the simplex $\sum_j \xi_{\ell,j} = 1$ is

$$f(\boldsymbol{\xi}) = \Gamma\left(\frac{n - m_0}{m_1}\right)\left[\prod_j \Gamma\left(\delta_{\ell,j}\frac{n - m_0}{m_1}\right)\right]^{-1}\prod_j \xi_{\ell,j}^{\delta_{\ell,j}((n-m_0)/m_1)}$$

so that the $\xi_{\ell,j}\hat{V}_\ell$ are independent gamma random variables with $\alpha = \delta_{\ell,j}((n - m_0)/2m_1)$.

The variance terms for $\mathcal{P}_2^*$ are constructed via

$$\hat{V}_j^* = \frac{m_0}{m_1}[\xi_{\ell,j}\hat{V}_\ell + \xi_{\ell+1,j}\hat{V}_{\ell+1}] \qquad \text{for } \frac{2\ell - 1}{2m_1} < \frac{j}{m_0} \le \frac{2\ell + 1}{2m_1},$$

which is a sum of gamma random variables. The weighting parameters are

$$\delta_{\ell,j} = \frac{m_1}{m_0}\left(\frac{2\ell + 1}{2} - \frac{2j - 1}{2}\left[\frac{m_1}{m_0}\right]\right) = \frac{m_1}{m_0}\zeta_{\ell,j},$$

$$\delta_{\ell+1,j} = \frac{m_1}{m_0}\left(\frac{2j - 1}{2}\left[\frac{m_1}{m_0}\right] - \frac{2\ell - 1}{2}\right) = \frac{m_1}{m_0}\zeta_{\ell+1,j}.$$

Thus $\frac{m_0}{m_1}\xi_{\ell,j}\hat{V}_\ell \sim \Gamma(\frac{n - m_0}{2m_0}\zeta_{\ell,j}, \frac{2m_0\sigma_\ell^2}{n - m_0})$.

On the edges of the interval, for $j/m_0 \le 1/2m_1$ and $j/m_0 > 1 - 1/2m_1$, the weights are simply $\delta_{1,j} = \delta_{m_1,j} = m_1/m_0$. There is no smoothing on the edges of the interval.

This is a somewhat involved transformation, and the divergence between the generated observations and the $\mathcal{P}_2^*$ observations is

$$(7.2) \qquad \sum_{j=1}^{m_0}\mathbf{D}(X_{0,j}, X_{0,j}^*) + \sum_{j=1}^{m_0}\mathbf{D}\left(\frac{m_0}{m_1}[\xi_{\ell,j}\hat{V}_\ell + \xi_{\ell+1,j}\hat{V}_{\ell+1}], V_j^*\right).$$

The first term can be bounded by noting that the logarithm of the variance function has a derivative bounded by $M$, and thus

$$(7.3) \qquad |\log \bar{\sigma}_j^2 - \log \sigma_\ell^2| = \zeta_{\ell+1,j}|\log \sigma_{\ell+1}^2 - \log \sigma_\ell^2| \le \frac{M}{m_1},$$

which, along with (10.1), implies

$$(7.4) \quad \mathbf{D}(X_{0,j}, X_{0,j}^*) = \frac{1}{2}\left[\frac{\bar{\sigma}_j^2}{\sigma_\ell^2} - 1 - \log\frac{\bar{\sigma}_j^2}{\sigma_\ell^2}\right] \le (\log \bar{\sigma}_j^2 - \log \sigma_\ell^2)^2 \le M^2 m_1^{-2}.$$

Thus the total error from the first term in (7.2) is less than $m_0 M^2 m_1^{-2}$. The contribution to the error from the edges is zero because $\bar{\sigma}_j^2 = \sigma_\ell^2$ for $j/m_0 \leq 1/2m_1$ or $j/m_0 \geq 1 - 1/2m_1$.

For the second term in (7.2), $\hat{V}_j^*$ is a sum of gammas with different scale terms. We need a general lemma on the distribution of gamma sums.

LEMMA 5. *For independent* $X_1, \ldots, X_m$ *random variables with* $X_i \sim \Gamma(\delta_i n, \sigma_i^2/n)$ *where the* $\delta_i > 0$ *and* $\sum_i \delta_i = 1$, *the distribution of the sum of the* $X_i$'s *is approximately gamma,*

$$\sum_i X_i \approx Y \sim \Gamma\left(n, \frac{\bar{\sigma}^2}{n}\right) \qquad \text{where } \bar{\sigma}^2 = \prod_{i=1}^m \sigma_i^{2\delta_i}.$$

*For* $r_i = \log \sigma_i^2 - \log \bar{\sigma}^2$, *the divergence between the distributions is bounded as* $r_i \to 0$ *by*

$$\mathbf{D}\left(\sum X_i, Y\right) \leq \sum_{i=1}^m \frac{n \delta_i r_i^4}{8} + \frac{r_i^2}{4} + O(n|r_i|^5 + |r_i|^3 + r_i^2 n^{-1}).$$

The proof of this lemma is in Section 11. For $r_\ell = \log(\sigma_\ell^2/\bar{\sigma}_j^2)$, Lemma 5 implies

$$\mathbf{D}(\hat{V}_j^*, V_j^*) \leq \left(\frac{n - m_0}{2m_0}\right)(\zeta_{\ell,j} r_\ell^4 + \zeta_{\ell+1,j} r_{\ell+1}^4) + \frac{r_\ell^2 + r_{\ell+1}^2}{4} + \cdots.$$

Using (7.3) once again, the bounds on the divergences are

$$\mathbf{D}(\hat{V}_j^*, V_j^*) \leq \left(\frac{n - m_0}{2m_0}\right)\frac{M^4}{m_1^4} + \frac{M^2}{2m_1^2} + O(nm_1^{-5} + m_1^{-3}).$$

On the edges the relationship is exact, $V_j^* = \xi_{\ell,j} \hat{V}_\ell$ for $j/m_0 \leq 1/2m_1$ and $j/m_0 > 1 - 1/2m_1$ where $\ell = 1$ and $\ell = m_1$, respectively. Therefore, the second term in (7.2) is bounded,

$$(7.5) \quad \begin{aligned} \sum_{j=1}^{m_0} \mathbf{D}(\hat{V}_j^*, V_j^*) &\leq \sum_{j=m_0/2m_1}^{m_0 - m_0/2m_1} \left(\frac{n - m_0}{m_0}\right)\frac{M^4}{m_1^4} + \frac{M^2}{m_1^2} \\ &\leq M^4 n m_1^{-4} + M^2 m_0 m_1^{-2}. \end{aligned}$$

Putting the two divergence bounds from (7.4) and (7.5) into (7.2) gives

$$(7.6) \quad \delta(\mathcal{P}_1^*, \mathcal{P}_2^*) \leq 2M m_0^{1/2} m_1^{-1} + M^2 n^{1/2} m_1^{-2}.$$

7.2.2. *The transformation of the $V_j^*$.* To reproduce the $\mathcal{P}_1^*$ random variables from the observation of $\mathcal{P}_2^*$, the chi-squared random variables $V_j^*$ are added together to generate approximately the $\hat{V}_\ell$.

As before, the transformation leaves the $X_{0,j}$ alone and the difference in the distributions of the $X_{0,j}^*$ and $X_{0,j}$ contributes a $O(m_0/m_1^2)$ term to the bound.

The $\hat{V}_\ell$ can be approximated by the sum of $V_j^*$ in the set $\mathcal{J}_\ell = \{j : \frac{\ell-1}{m_1} < \frac{j}{m_0} \le \frac{\ell}{m_1}\}$. By Lemma 5, the sum of the $\frac{m_1}{m_0} V_j^*$ will be approximately a gamma with expectation

$$\log \sigma_\ell^{*2} = \sum_{j=m_0(\ell-1)/m_1+1}^{m_0(\ell-1/2)/m_1} (\zeta_{\ell-1,j} \log \sigma_{\ell-1}^2 + \zeta_{\ell,j} \log \sigma_\ell^2)$$

$$+ \sum_{j=m_0(\ell-1/2)/m_1+1}^{m_0\ell/m_1} (\zeta_{\ell,j} \log \sigma_\ell^2 + \zeta_{\ell+1,j} \log \sigma_{\ell+1}^2),$$

which equals $\frac{1}{8} \log \sigma_{\ell-1}^2 + \frac{3}{4} \log \sigma_\ell^2 + \frac{1}{8} \log \sigma_{\ell+1}^2$. The correction at the edges implies that $\log \sigma_1^{*2} = \frac{7}{8} \log \sigma_1^2 + \frac{1}{8} \log \sigma_2^2$ and $\log \sigma_{m_1}^{*2} = \frac{1}{8} \log \sigma_{m_0-1}^2 + \frac{7}{8} \log \sigma_{m_0}^2$.

To the error bounded in Lemma 5, we will have to add the error from the difference between $\sigma_\ell^{*2}$ and $\sigma_\ell^2$. Let $g_\sigma(x)$ be the density of a gamma distribution with $\alpha = (n-m_0)/2m_1$ and $\beta = 2m_1\sigma^2/(n-m_0)$. The divergence can then be written as

$$
\begin{aligned}
(7.7) \quad \mathbf{D}\left(\frac{m_1}{m_0} \sum_{j \in \mathcal{J}_\ell} V_j^*, \hat{V}_\ell\right) = {} & \mathbf{D}\left(\frac{m_1}{m_0} \sum_{j \in \mathcal{J}_\ell} V_j^*, \check{V}_\ell\right) \\
& + \mathbb{E} \log \frac{g_{\sigma_\ell^*}(m_1/m_0 \sum_{j \in \mathcal{J}_\ell} V_j^*)}{g_{\sigma_\ell}(m_1/m_0 \sum_{j \in \mathcal{J}_\ell} V_j^*)},
\end{aligned}
$$

where $\check{V}_\ell$ is the gamma random variable with density $g_{\sigma_\ell^*}$,

Let $r_j = \log \bar{\sigma}_j^2 / \sigma_\ell^2$. Then by Lemma 5,

$$\mathbf{D}\left(\frac{m_1}{m_0} \sum_{j \in \mathcal{J}_\ell} V_j^*, \check{V}_\ell\right) \le \sum_{j \in \mathcal{J}_\ell} \left[\left(\frac{n-m_0}{2m_0}\right) \frac{r_j^4}{8} + \frac{r_j^2}{4} + \cdots\right],$$

where the $r_j$ are bounded for $(\ell-1)/m_1 < j/m_0 < (2\ell-1)/2m_1$ using (7.3),

$$|r_j| = \left|\frac{1}{8} \log \frac{\sigma_\ell^2}{\sigma_{\ell-1}^2} + \left(\zeta_{\ell,j} - \frac{7}{8}\right) \log \frac{\sigma_\ell^2}{\sigma_{\ell+1}^2}\right| \le \frac{M}{m_1}.$$

For $(2\ell-1)/2m_1 < j/m_0 \le \ell/m_1$, there is an analogous calculation so that $|r_j| < M/m_1$ for every $j$. For $j/m_0 < 1/2m_1$, the $|r_j| = |\log \sigma_1^2/\sigma_2^2|/8 \le M/m_1$. There is an analogous bound on the other end of the interval as well.

Thus,

$$(7.8) \quad \sum_\ell \mathbf{D}\left(\frac{m_1}{m_0}\sum_{j\in\mathcal{J}_\ell} V_j^*, \check{V}_\ell\right) \leq \frac{M^4(n-m_0)}{16m_1^4} + \frac{M^2 m_0}{4m_1^2} + O(m_0 m_1^{-3}).$$

The second term in (7.7) is the expected value of

$$\log\frac{g_{\sigma_\ell^*}(X)}{g_{\sigma_\ell}(X)} = \frac{n-m_0}{2m_1}\left[\log\frac{\sigma_\ell^2}{\sigma_\ell^{*2}} + X\left(\frac{1}{\sigma_\ell^2} - \frac{1}{\sigma_\ell^{*2}}\right)\right].$$

The expectation is to be taken over the distribution of the average of the $V_j^*$,

$$\mathbb{E}\frac{m_1}{m_0}\sum_{j\in\mathcal{J}_\ell} V_j^* = \frac{m_1}{m_0}\sum_{j\in\mathcal{J}_\ell}\bar\sigma_j^2 = \frac{m_1}{m_0}\sum_{j\in\mathcal{J}_\ell}\sigma_\ell^{\zeta_{\ell,j}}\sigma_{\ell+1}^{1-\zeta_{\ell,j}}.$$

Section 13 does the necessary calculation to bound the contribution from this expectation. From (13.2) and (13.3),

$$\sum_{\ell=1}^{m_1}\mathbb{E}\log\frac{g_{\sigma_\ell^*}(X)}{g_{\sigma_\ell}(X)} \leq M^2 n m_1^{-3} + \frac{M^2}{m_1}\sum_{\ell=2}^{m_1-1} n m_1^{-2\alpha_1}.$$

Putting this together with the bound in (7.8),

$$\delta(\mathcal{P}_2^*, \mathcal{P}_1^*) \leq M n^{1/2} m_1^{-\alpha_1} + M m_0^{1/2} m_1^{-1} + M n^{1/2} m_1^{-3/2}.$$

Therefore, in light of the analogous result in (7.6),

$$(7.9) \quad \Delta(\mathcal{P}_2^*, \mathcal{P}_1^*) \leq M n^{1/2} m_1^{-\alpha_1} + 2M m_0^{1/2} m_1^{-1} + M n^{1/2} m_1^{-3/2}.$$

7.3. *Approximating the nonparametric regression.* The last step is to show that $\mathcal{P}_2^*$ is equivalent to the $n$ independent normal observations from the original $\tilde{\mathcal{P}}$ experiment.

The $V_j^*$ observations are sufficient statistics for $(n-m_0)/m_0$ independent normals with means 0 and variances $\bar\sigma_j^2$. These normals will be used in place of all the wavelet coefficients $X_{i,j}$. These wavelet coefficients are combined with the $X_{0,j}^*$, and using Mallat's algorithm, $n$ normal observations are produced, $Y_i^* \sim \mathcal{N}(\sqrt{m_0}\theta_{k_0,j}, \bar\sigma_j^2)$ where $(j-1)/m_0 < i/n \leq j/m_0$.

The error made by this approximation is

$$\mathbf{D}(Y_i^*, Y_i) = \frac{(\sqrt{m_0}\theta_{k_0,j} - f(i/n))^2}{2\bar\sigma_j^2} + \frac{1}{2}\left[\frac{\sigma^2(i/n)}{\bar\sigma_j^2} - 1 - \log\left(\frac{\sigma^2(i/n)}{\bar\sigma_j^2}\right)\right].$$

The mean functions are bounded much as with the constant variance case,

$$(7.10) \quad \sum_i \frac{(f(i/n) - \sqrt{m_0}\theta_{k_0,j})^2}{\bar\sigma_j^2} \leq e^M \frac{n}{m_0^{2\alpha}}\left(\sum_{i\geq k_0} 2^{i(\alpha+1/2)}\sup_j\left|\frac{\theta_{i,j^*}}{\bar\sigma}\right|\right)^2$$

$$\leq e^M \frac{n}{m_0^{2\alpha}}\gamma_{k_0}^2,$$

because the partial sums of the $b(\alpha, \infty, 1)$ norm are uniformly bounded for $f \in \tilde{\mathcal{F}}_{\bar{\sigma}}(\alpha, \gamma_k)$ and using the smoothness of the variances for $\bar{\sigma}^2 \leq e^M \bar{\sigma}_j^2$.

The second part of the divergence is

$$\sum_{i=1}^{n} \frac{1}{2} \left[ \frac{\sigma^2(i/n)}{\bar{\sigma}_j^2} - 1 - \log\left(\frac{\sigma^2(i/n)}{\bar{\sigma}_j^2}\right) \right] \leq \sum_{i=1}^{n} (\log \sigma^2(i/n) - \log \bar{\sigma}_j^2)^2.$$

To bound this quantity requires taking advantage of all the smoothness in the variance functions. To simplify things a bit, we write $\tau(t) = \log \sigma^2(t)$ and $t_j^* = (2j-1)/(2m_0)$. The smoothness condition on $\tau$ implies that $\tau(t) = \tau(t_j^*) + (t - t_j^*)\tau'(t_j^*) + E$ where the error term is $|E| \leq M|t - t_j^*|^\alpha$. By the definition of $\bar{\sigma}_j^2$,

$$\log \bar{\sigma}_j^2 = m_1 \int_{(\ell-1)/m_1}^{\ell/m_1} \zeta_{\ell,j} \tau(t) + \zeta_{\ell+1,j} \tau(t + 1/m_1) \, dt$$

$$= \tau(t_j^*)$$
$$\quad + m_1 \tau'(t_j^*) \int_{(\ell-1)/m_1}^{\ell/m_1} [\zeta_{\ell,j}(t - t_j^*) + \zeta_{\ell+1,j}(t - t_j^* + 1/m_1)] \, dt + E_1$$

$$= \tau(t_j^*)$$
$$\quad + \tau'(t_j^*) \left[ \zeta_{\ell,j}\left(\frac{2\ell-1}{2m_1} - \frac{2j-1}{2m_0}\right) + \zeta_{\ell+1,j}\left(\frac{2\ell+1}{2m_1} - \frac{2j-1}{2m_0}\right) \right] + E_1.$$

The error is an average over errors in the expansion and so $|E_1| \leq M m_1^{-\alpha_1}$. Plugging in $\zeta_{\ell,j}$ and $\zeta_{\ell+1,j}$ according to their definitions,

$$\zeta_{\ell,j}\left(\frac{2\ell-1}{2m_1} - \frac{2j-1}{2m_0}\right) + \zeta_{\ell+1,j}\left(\frac{2\ell+1}{2m_1} - \frac{2j-1}{2m_0}\right)$$

$$= -\frac{\zeta_{\ell,j}\zeta_{\ell+1,j}}{m_1} + \frac{\zeta_{\ell+1,j}\zeta_{\ell,j}}{m_1} = 0.$$

Thus, $\log \bar{\sigma}_j^2 = \tau(t_j^*) + E_1$, and, from the bound on the derivative, $|\tau(\frac{i}{n}) - \tau(t_j^*)| \leq M m_0^{-1}$ implies that the bound is

$$(7.11) \quad \left[ \sum_{i=1}^{n} (\log \sigma^2(i/n) - \log \bar{\sigma}_j^2)^2 \right]^{1/2} \leq M n^{1/2} m_1^{-\alpha_1} + M n^{1/2} m_0^{-1}.$$

Combining (7.11) and (7.10) implies

$$(7.12) \quad \Delta(\mathcal{P}_2^*, \breve{\mathcal{P}}) \leq e^{M/2} n^{1/2} m_0^{-\alpha} \gamma_{k_0} + M n^{1/2} m_1^{-\alpha_1} + M n^{1/2} m_0^{-1}.$$

It is not necessary to do a specific calculation to bound $\delta(\mathcal{P}, \mathcal{P}_2^*)$ because Mallat's algorithm is invertible, and the deficiency distance between the distributions after applying the inverse to both distributions can only be smaller.

Putting together (7.1), (7.9) and (7.12),

$$\Delta(\bar{\mathcal{P}}_k, \breve{\mathcal{P}}) \leq 2e^{M/2}n^{1/2}m_0^{-\alpha}\gamma_{k_0}$$
$$+ 2Mn^{1/2}m_1^{-\alpha_1} + 2Mm_0^{1/2}m_1^{-1} + Mn^{1/2}m_0^{-1} + Mn^{1/2}m_1^{-3/2},$$

which proves Lemma 3.

## 8. The variance process.

The final piece of the proof of Theorem 2 is to show that the variance observations in the simplified experiment can be transformed into a continuous Gaussian process. The experiment $\tilde{\mathcal{Q}}_n$ observes $m_1$ independent variance components $V_\ell$ and a countable sequence of normal coefficients. Mimicking the results of [10] and [24], we can construct an independent Gaussian process $V(t)$ from the $\chi^2$ observations, $dV(t) = \log\sigma^2(t)\,dt + \sqrt{2}n^{-1/2}\,dW_2(t)$. The construction follows by taking the logarithm of the $V_\ell$ and then using them to approximate the increments of $V(t)$.

Taking the logarithm of the $V_\ell$ generates an intermediate experiment $\mathcal{Q}^*$ that observes $Z_\ell \sim \mathcal{N}(\log\sigma_\ell^2, \frac{2m_1}{n})$ all independent and the Gaussian process, conditional on the $Z_\ell$,

$$dY^*(t) = f(t)\,dt + e^{Z_\ell/2}n^{-1/2}\,dW(t) \qquad \text{for } \frac{\ell-1}{m_1} \leq t \leq \frac{\ell}{m_1}.$$

The divergence between the distributions is

$$\mathbf{D}((\log V, Y), (Z, Y^*)) = \sum_{\ell=1}^{m_1}\mathbf{D}(\log V_\ell, Z_\ell) + \mathbf{D}(Y, Y^*|\log V).$$

The first term is bounded as in Section 10.1 by $m_1^2/n$ and for the second the divergence term $\mathbf{D}(Y, Y^*|V) = 0$ because the conditional distributions are the same.

In the other direction, the $V_\ell$ are approximated by $\exp[Z_\ell]$. The divergence bounds are the same (these transformations are one-to-one and increasing), thus

$$(8.1) \qquad\qquad \Delta(\tilde{\mathcal{Q}}, \mathcal{Q}^*) \leq m_1 n^{-1/2}.$$

The scaled increments $m_1[V(\ell/m_1) - V([\ell-1]/m_1)]$ from $\breve{\mathcal{Q}}$ have the same distribution as the $Z_\ell$ from $\mathcal{Q}^*$. Thus, $\delta(\breve{\mathcal{Q}}, \mathcal{Q}^*) = 0$.

To bound $\delta(\mathcal{Q}^*, \breve{\mathcal{Q}})$ we need to construct the entire $V(t)$ process via a smoothing operation on the $Z_\ell$ that is described in detail in [3]. Our argument follows this reference very closely, so only an outline of the steps will be given.

The transformation uses triangular interpolating kernels,

$$K_\ell(x) = m_1 - m_1^2 \left| x - \frac{2\ell - 1}{2m_1} \right| \quad \text{for } \frac{2\ell - 3}{2m_1} \le x \le \frac{2\ell + 1}{2m_1}$$

with the appropriate reflections at the boundaries. The variable $V(t)$ is constructed from the $Z_\ell$,

$$dV^*(t) = \sum_{\ell=1}^{m_1} Z_\ell K_\ell(t) + \sqrt{2}\, n^{-1/2} \sum_{\ell=1}^{m_1} \frac{1}{\sqrt{m_1}}\, dB_\ell(t),$$

where the $B_\ell(t)$ are independent $K_\ell$-Brownian bridges. This Gaussian process is (see [3])

$$dV^*(t) = \hat{\tau}(t)\, dt + \sqrt{2}\, n^{-1/2}\, dW_2(t),$$

where $W_2(t)$ is a standard Brownian motion and $\hat{\tau}$ is a piecewise-linear function where $\hat{\tau}(t_\ell^*) = \log \sigma_\ell^2$ for $t_\ell^* = \frac{2\ell - 1}{2m_1}$ and

$$(8.2) \qquad \hat{\tau}(t) = m_1(t_{\ell+1}^* - t)\hat{\tau}(t_\ell^*) + m_1(t - t_\ell^*)\hat{\tau}(t_{\ell+1}^*)$$

for $t_\ell^* \le t \le t_{\ell+1}^*$. For $t < \frac{1}{2m_1}$ or $t > 1 - \frac{1}{2m_1}$, the function $\hat{\tau}$ is just a constant equal to $\hat{\tau}(\frac{1}{2m_1})$ or $\hat{\tau}(1 - \frac{1}{2m_1})$, respectively.

The smoothness condition on the $\tau$ functions implies that

$$(8.3) \qquad \hat{\tau}(t_\ell^*) = m_1 \int_{(\ell-1)/m_1}^{\ell/m_1} \tau(x)\, dx = \tau(t) + m_1(t_\ell^* - t)\tau'(t) + \bar{E}_1.$$

Likewise, $\hat{\tau}(t_{\ell+1}^*) = \tau(t) + m_1(t_{\ell+1}^* - t)\tau'(t) + \bar{E}_2$. Therefore, plugging (8.3) into (8.2) yields $|\hat{\tau}(t) - \tau(t)| \le 2Mm_1^{-\alpha_1}$ for $t$ in the interior of the interval. Unfortunately, at the boundaries the error is of the order $m_1^{-1}$. Thus, the $L_2$ distance between $\hat{\tau}$ and $\tau$ is bounded by

$$(8.4) \qquad \|\hat{\tau} - \tau\|_2^2 \le 4M^2 m_1^{-2\alpha_1} + M^2 m_1^{-3}.$$

The total-variation distance between the distributions of $V^*(t)$ and $V(t)$ is of the order of this distance divided by the variance. Therefore, from (8.1) and (8.4),

$$(8.5) \qquad \Delta(\check{Q}, \tilde{Q}) \le m_1 n^{-1/2} + 2Mn^{1/2} m_1^{-\alpha_1} + Mn^{1/2} m_1^{-3/2}.$$

This proves Lemma 4.

## 9. Divergence bounds for gamma distributions.

LEMMA 6. *The $K$–$L$ divergence between $\mathbb{P}_1 = \Gamma(\alpha_1, \beta_1)$ and $\mathbb{P}_2 = \Gamma(\alpha_2, \beta_2)$ is*

$$(9.1) \qquad \mathbf{D}(\mathbb{P}_1, \mathbb{P}_2) \le \frac{(\alpha_1 - \alpha_2)^2}{2\alpha_1^2} + O\!\left(\frac{(\alpha_1 - \alpha_2)^2}{\alpha_1^3}\right)$$

*when the means are the same ($\alpha_1 \beta_1 = \alpha_2 \beta_2$).*

The divergence between a pair of gamma distributions is

$$(9.2) \quad \mathbf{D}(\mathbb{P}_1, \mathbb{P}_2) = \frac{\alpha_1(\beta_1 - \beta_2)}{\beta_2} + \alpha_2 \log\left(\frac{\beta_2}{\beta_1}\right)$$
$$+ \log\left(\frac{\Gamma(\alpha_2)}{\Gamma(\alpha_1)}\right) + (\alpha_1 - \alpha_2)\mathbb{P}_1 \log\left[\frac{x}{\beta_1}\right].$$

We can rewrite (9.2) just in terms of the $\alpha$'s using the substitution $\frac{\beta_1}{\beta_2} = \frac{\alpha_2}{\alpha_1}$,

$$(9.3) \quad \mathbf{D}(\mathbb{P}_1, \mathbb{P}_2) = (\alpha_2 - \alpha_1)$$
$$+ \alpha_2 \log\left(\frac{\alpha_1}{\alpha_2}\right) + \log\left(\frac{\Gamma(\alpha_2)}{\Gamma(\alpha_1)}\right) + (\alpha_1 - \alpha_2)\mathbb{P}_1 \log\left[\frac{x}{\beta_1}\right].$$

Let $\delta = \alpha_2 - \alpha_1$. As in the proof of Lemma 7, the last two terms can be bounded using the integral remainder of a Taylor series,

$$\log\left(\frac{\Gamma(\alpha_1 + \delta)}{\Gamma(\alpha_1)}\right) - \delta\mathbb{P}_1 \log\left[\frac{x}{\beta_1}\right] = \int_0^\delta \psi'(\alpha_1 + t)(\delta - t)\,dt.$$

On the other hand, the first two terms in (9.3) have a similar Taylor series form,

$$(\alpha_2 - \alpha_1) + \alpha_2 \log\frac{\alpha_1}{\alpha_2} = \delta - (\alpha_1 + \delta)\log[1 + \delta/\alpha_1] = -\int_0^\delta \frac{\delta - t}{\alpha_1 + t}\,dt.$$

The classical expansion in (12.2) implies that the K–L divergence is

$$\mathbf{D}(\mathbb{P}_1, \mathbb{P}_2) = \int_0^\delta \psi'(\alpha_1 + t)(\delta - t)\,dt - \int_0^\delta \frac{\delta - t}{\alpha_1 + t}\,dt$$
$$= \int_0^\delta \frac{\delta - t}{2(\alpha_1 + t)^2} + (\delta - t)O(\alpha_1^{-3})\,dt,$$

which gives the bound asserted by the lemma.

9.1. *What if we increase the degrees of freedom?* In order to take into account the noncentrality of some $\chi^2$ distributions, we need a bound on the divergence between a $\Gamma(\alpha_1, \beta_1)$ and a $\mathbb{P}_2^* = \Gamma(\alpha_2 + \lambda, \beta_2)$ when $\alpha_1\beta_1 = \alpha_2\beta_2$ and $\lambda > 0$,

$$(9.4) \quad \mathbf{D}(\mathbb{P}_1, \mathbb{P}_2^*) = \mathbf{D}(\mathbb{P}_1, \mathbb{P}_2) + \lambda \log\left(\frac{\alpha_1}{\alpha_2}\right)$$
$$+ \log\left(\frac{\Gamma(\alpha_2 + \lambda)}{\Gamma(\alpha_2)}\right) - \lambda\mathbb{P}_1 \log\left[\frac{x}{\beta_1}\right].$$

From Lemma 7

$$(9.5) \quad \log\left(\frac{\Gamma(\alpha_2 + \lambda)}{\Gamma(\alpha_2)}\right) - \lambda\mathbb{P}_2 \log\left[\frac{x}{\beta_1}\right] = \lambda^2\left(\frac{1}{2\alpha_2} + O(\alpha_2^{-2})\right),$$

and the difference in taking the expectation with respect to $\mathbb{P}_1$ instead of $\mathbb{P}_2$ can be bounded using Jensen's inequality,

$$\lambda \log\left(\frac{\alpha_1}{\alpha_2}\right) + \lambda\left[\mathbb{P}_2 \log\left[\frac{x}{\beta_1}\right] - \mathbb{P}_1 \log\left[\frac{x}{\beta_1}\right]\right] \leq \lambda \log\frac{\alpha_1}{\alpha_1 - 1} \leq \frac{\lambda}{\alpha_1 - 1},$$

Therefore, substituting this last inequality, (9.5), and Lemma 6 into (9.4),

$$
\begin{aligned}
(9.6) \qquad \mathbf{D}(\mathbb{P}_1, \mathbb{P}_2^*) \leq & \frac{(\alpha_1 - \alpha_2)^2}{2\alpha_1^2} + \frac{\lambda^2}{2\alpha_2} + \frac{\lambda}{\alpha_1 - 1} \\
& + O(\lambda^2 \alpha_2^{-2} + (\alpha_1 - \alpha_2)^2 \alpha_1^{-3}).
\end{aligned}
$$

**10. Divergence between normal distributions.** If there are two normal distributions with means $\mu_1$ and $\mu_2$, and variances $\sigma_1^2$ and $\sigma_2^2$, respectively, then the divergence between them is

$$(10.1) \qquad \mathbf{D}(N_1, N_2) = \frac{1}{2}\left[\frac{\sigma_1^2}{\sigma_2^2} - 1 - \log\left(\frac{\sigma_1^2}{\sigma_2^2}\right)\right] + \frac{(\mu_1 - \mu_2)^2}{2\sigma_2^2}.$$

10.1. *The logarithm of the gamma distribution.* Suppose that $X$ has a $\Gamma(\alpha, 1)$ distribution and $W = \log(X)$. Let $f(w)$ be the density of $W$ which is approximately normal for large $\alpha$. Let $\phi_\alpha(z)$ be the density of a $Z \sim \mathcal{N}(\log\alpha, \alpha^{-1})$ distribution. The K–L divergence between these distributions can be bounded using Stirling's formula,

$$\mathbf{D}(Z, W) = \log\Gamma(\alpha) - \log(\sqrt{2\pi}) + \frac{1}{2}\log\alpha + \alpha e^{1/(2\alpha)} - \alpha\log\alpha - \frac{1}{2}$$

$$= \alpha\left(\exp\left[\frac{1}{2\alpha}\right] - 1\right) - \frac{1}{2} + \frac{\theta}{12\alpha} \leq \frac{1}{3\alpha} \qquad \text{for } \alpha \geq 1/2.$$

To extend this bound, notice the logarithm of a $\Gamma(\alpha, \beta)$ random variable is a shift of the distribution of $\log X$. Therefore, if $X \sim \Gamma(n/2, 2\sigma^2/n)$, then $\mathbf{D}(Z, \log(X)) \leq n^{-1}$ where $Z \sim \mathcal{N}(\log\sigma^2, \frac{2}{n})$. Compare this to [14], Lemma A.3.

**11. Proof of Lemma 5.** To bound the divergence between the sums, we define some similar random variables $X_i^* \sim \Gamma(\delta_i(1 + r_i)n, \frac{\bar{\sigma}^2}{n})$. The definition of $r_i$ implies that $\sum \delta_i r_i = 0$, and thus the distribution of the sum of the $X_i^*$'s is the same as the distribution of $Y$. It is necessary that $1 + r_i > 0$, but the bound is only interesting for small $r_i$ anyway.

We bound the divergence between the sums by the divergence between the joint distributions, $\mathbf{D}(\sum X_i, \sum X_i^*) \leq \mathbf{D}((X_1, \ldots, X_m), (X_1^*, \ldots, X_m^*)) = \sum_{i=1}^m \mathbf{D}(X_i, X_i^*)$.

The divergence can then be bounded using

$$\sum_{i=1}^{m} \mathbf{D}(X_i, X_i^*) = \sum_{i=1}^{m} \log\left[\frac{\Gamma(\delta_i(1+r_i)n)}{\Gamma(\delta_i n)}\right]$$

$$(11.1) \hspace{3cm} - \delta_i r_i n \mathbb{E} \log X_i + \mathbb{E}\left[nX_i\left(\frac{1}{\bar{\sigma}^2} - \frac{1}{\sigma_i^2}\right)\right].$$

From Lemma 7, the first two terms of (11.1) are

$$\log\left[\frac{\Gamma(\delta_i n + \delta_i r_i n)}{\Gamma(\delta_i n)}\right] - \delta_i r_i n \mathbb{E} \log X_i$$

$$(11.2) \hspace{1cm} = -\delta_i r_i n \log\left[\frac{\sigma_i^2}{n}\right] + \frac{n\delta_i r_i^2}{2}$$

$$+ \frac{r_i^2}{4} - \frac{n\delta_i r_i^3}{6} + \frac{n\delta_i r_i^4}{12} + O\left(\frac{\delta_i r_i^2}{n} + |r_i|^3 + \delta_i |r_i|^5\right).$$

In summing over $i$, the first term is $n\sum_{i=1}^{m} \delta_i r_i (\log \sigma_i^2 - \log n) = n\sum_{i=1}^{m} \delta_i r_i^2$ because $\sum \delta_i r_i \log \bar{\sigma}^2 = \sum \delta_i r_i \log n = 0$. Summing over all the terms in (11.2) yields a bound on the contribution from the first two terms of (11.1),

$$(11.3) \quad \sum_i \left[-\frac{n\delta_i r_i^2}{2} + \frac{r_i^2}{4} - \frac{n\delta_i r_i^3}{6} + \frac{n\delta_i r_i^4}{12} + O(n\delta_i |r_i|^5 + |r_i|^3 + \delta_i r_i^2 n^{-1})\right].$$

The last term in (11.1) is

$$\mathbb{E}\left[nX\left(\frac{1}{\bar{\sigma}^2} - \frac{1}{\sigma_i^2}\right)\right] = n\delta_i(e^{r_i} - 1)$$

$$(11.4) \hspace{2cm} = n\delta_i r_i + \frac{n\delta_i r_i^2}{2} + \frac{n\delta_i r_i^3}{6} + \frac{n\delta_i r_i^4}{24} + O(n\delta_i |r_i|^5).$$

By summing (11.4) over $i$ and then adding it to (11.3), we bound the divergence by

$$\mathbf{D}\left(\sum X_i, \sum X_i^*\right) \le \sum_{i=1}^{m} \frac{n\delta_i r_i^4}{8} + \frac{r_i^2}{4} + O(n\delta_i |r_i|^5 + |r_i|^3 + \delta_i r_i^2 n^{-1}).$$

## 12. Digamma bound.

LEMMA 7. *For* $X \sim \Gamma(\alpha, 1)$,

$$\log\left(\frac{\Gamma(\alpha + \delta)}{\Gamma(\alpha)}\right) - \delta \mathbb{E} \log X = \frac{\delta^2}{2}\left(\frac{1}{\alpha} + \frac{1}{2\alpha^2} + O(\alpha^{-3})\right) - \frac{\delta^3}{6}\left(\frac{1}{\alpha^2} + O(\alpha^{-3})\right)$$

$$+ \frac{\delta^4}{24}\left(\frac{2}{\alpha^3} + O(\alpha^{-4})\right) + O(\delta^5 \alpha^{-4})$$

*as* $\alpha \to \infty$ *and* $\delta/\alpha \to 0$.

This lemma is essentially a result on the properties of gamma and "polygamma" functions,

$$\Gamma'(\alpha) = \int_0^\infty \frac{d}{d\alpha}\, x^{\alpha-1} e^{-x}\, dx = \int_0^\infty (\log x) x^{\alpha-1} e^{-x}\, dx = \Gamma(\alpha)\mathbb{E}\log X.$$

Thus,

$$\log\!\left(\frac{\Gamma(\alpha+\delta)}{\Gamma(\alpha)}\right) - \delta\mathbb{E}\log X = \log\!\left(\frac{\Gamma(\alpha+\delta)}{\Gamma(\alpha)}\right) - \delta\frac{\Gamma'(\alpha)}{\Gamma(\alpha)}.$$

The $(n+1)$th derivative of the logarithm of the gamma function is the polygamma function $\psi^{(n)}(\alpha)$. Thus the expression can be seen as a Taylor expansion of $\log\Gamma(k)$ around $k=\alpha$,

$$
\begin{aligned}
(12.1)\qquad \log\!\left(\frac{\Gamma(\alpha+\delta)}{\Gamma(\alpha)}\right) &= \delta\psi^{(0)}(\alpha) + \frac{\delta^2}{2}\psi^{(1)}(\alpha) \\
&\quad + \frac{\delta^3}{6}\psi^{(2)}(\alpha) + \frac{\delta^4}{24}\psi^{(3)}(\alpha) + O(\delta^5\psi^{(4)}(\alpha)).
\end{aligned}
$$

There is a classical result,

$$(12.2)\qquad \psi^{(1)}(\alpha) = \frac{1}{\alpha} + \frac{1}{2\alpha^2} + O(\alpha^{-3}),$$

and plugging the derivatives of this function into (12.1) gives the desired expansion.

**13. Bound.** In Section 7.2.2, we need a bound on the quantity

$$S_\ell = \mathbb{E}\log\frac{g_{\sigma_\ell^*}(X)}{g_{\sigma_\ell}(X)} = \frac{n-m_0}{2m_0}\sum_{j\in\mathcal{J}_\ell}\log\!\left[\frac{\sigma_\ell^2}{\sigma_\ell^{*2}}\right] + \frac{\bar{\sigma}_j^2}{\sigma_\ell^2} - \frac{\bar{\sigma}_j^2}{\sigma_\ell^{*2}}.$$

To straighten out this sum we need to separate it into two sets of $j$'s,

$$\mathcal{J}_\ell^1 = \left\{j : \frac{\ell-1}{m_1} < \frac{j}{m_0} \le \frac{2\ell-1}{2m_1}\right\} \quad\text{and}\quad \mathcal{J}_\ell^2 = \left\{j : \frac{2\ell-1}{2m_1} < \frac{j}{m_0} \le \frac{\ell}{m_1}\right\},$$

so that $\mathcal{J}_\ell^1 \cup \mathcal{J}_\ell^2 = \mathcal{J}_\ell$. For $j\in\mathcal{J}_\ell^1$ the variances $\sigma_j^{2*}$ are $\sigma_\ell^{2\zeta_{\ell,j}}\sigma_{\ell-1}^{2\zeta_{\ell-1,j}}$, and for $j\in\mathcal{J}_\ell^2$ the variances $\sigma_\ell^{2*}$ are $\sigma_\ell^{2\zeta_{\ell,j}}\sigma_{\ell+1}^{2\zeta_{\ell+1,j}}$. Thus the sum can be written as

$$
\begin{aligned}
S_\ell = \frac{n-m_0}{2m_0}\sum_{j\in\mathcal{J}_\ell^1}\log\!\left[\frac{\sigma_\ell^2}{\sigma_{\ell-1}^{1/4}\sigma_\ell^{3/2}\sigma_{\ell+1}^{1/4}}\right] &+ \left(\frac{\sigma_{\ell-1}^2}{\sigma_\ell^2}\right)^{1-\zeta_{\ell,j}} - \left[\frac{\sigma_\ell^{2\zeta_{\ell,j}}\sigma_{\ell-1}^{2(1-\zeta_{\ell,j})}}{\sigma_{\ell-1}^{1/4}\sigma_\ell^{3/2}\sigma_{\ell+1}^{1/4}}\right] \\
+ \frac{n-m_0}{2m_0}\sum_{j\in\mathcal{J}_\ell^2}\log\!\left[\frac{\sigma_\ell^2}{\sigma_{\ell-1}^{1/4}\sigma_\ell^{3/2}\sigma_{\ell+1}^{1/4}}\right] &+ \left(\frac{\sigma_{\ell+1}^2}{\sigma_\ell^2}\right)^{1-\zeta_{\ell,j}} - \left[\frac{\sigma_\ell^{2\zeta_{\ell,j}}\sigma_{\ell+1}^{2(1-\zeta_{\ell,j})}}{\sigma_{\ell-1}^{1/4}\sigma_\ell^{3/2}\sigma_{\ell+1}^{1/4}}\right].
\end{aligned}
$$

Note that at the edges, where the $\ell = 1$ or $\ell = m_1$, this equation is still true if we define $\sigma_0^2 = \sigma_1^2$ and $\sigma_{m_1+1}^2 = \sigma_{m_1}^2$.

Setting $r_\ell = \log \sigma_{\ell+1}^2 - \log \sigma_\ell^2$ and $r_{\ell-1} = \log \sigma_\ell^2 - \log \sigma_{\ell-1}^2$, we can write $S_\ell$ as

$$S_\ell = \frac{n-m_0}{2m_0} \sum_{j \in \mathcal{J}_\ell^1} \frac{(r_{\ell-1} - r_\ell)}{8} + \exp[-r_{\ell-1}(1 - \zeta_{\ell,j})]\left[1 - \exp\left(\frac{1}{8}(r_{\ell-1} - r_\ell)\right)\right]$$

$$+ \frac{n-m_0}{2m_0} \sum_{j \in \mathcal{J}_\ell^2} \frac{(r_{\ell-1} - r_\ell)}{8} + \exp[r_\ell(1 - \zeta_{\ell,j})]\left[1 - \exp\left(\frac{1}{8}(r_{\ell-1} - r_\ell)\right)\right]$$

$$\leq \frac{n-m_0}{2m_0}\left[\sum_{j \in \mathcal{J}_\ell^1} \left(\frac{r_{\ell-1} - r_\ell}{8}\right)[1 - e^{-r_{\ell-1}(1 - \zeta_{\ell,j})}]\right.$$

$$+ \left.\sum_{j \in \mathcal{J}_\ell^2} \left(\frac{r_{\ell-1} - r_\ell}{8}\right)[1 - e^{r_\ell(1 - \zeta_{\ell,j})}]\right]$$

$$\leq \frac{n-m_0}{2m_0}\left(\frac{r_{\ell-1} - r_\ell}{8}\right)\left[\sum_{j \in \mathcal{J}_\ell^1} r_{\ell-1}(1 - \zeta_{\ell,j}) - \sum_{j \in \mathcal{J}_\ell^2} r_\ell(1 - \zeta_{\ell,j})\right].$$

By the definition of the weights $\zeta_{\ell,j}$, each one is between 0 and 1 and the sums of them are $\sum_{j \in \mathcal{J}_\ell^1}(1 - \zeta_{\ell,j}) = \sum_{j \in \mathcal{J}_\ell^2}(1 - \zeta_{\ell,j}) = \frac{3}{8}\left(\frac{m_0}{m_1}\right)$. Therefore,

$$(13.1) \qquad\qquad S_\ell \leq \frac{3(n-m_0)}{128 m_1}(r_{\ell-1} - r_\ell)^2.$$

The definition of the function class says that the function $\tau(t) = \log(\sigma^2(t))$ is smooth in the sense that $\tau(t+\delta) = \tau(t) + \delta\tau'(t) + E$ where $|E| \leq M\delta^{\alpha_1}$. To use this notice that, for $\ell = 1, \ldots, m_1 - 1$, expanding the functions around $(2\ell-1)/2m_1$,

$$r_\ell - r_{\ell-1} = m_1 \int_{-1/2m_1}^{1/2m_1} \tau\left(\frac{2\ell+1}{2m_1} + t\right) - 2\tau\left(\frac{2\ell-1}{2m_1} + t\right) + \tau\left(\frac{2\ell-3}{2m_1} + t\right) dt$$

$$= m_1 \int_{-1/2m_1}^{1/2m_1} E_1 - 2E_2 + E_3 \, dt,$$

and the average value of the errors is less than $m_1 \int_{-1/2m_1}^{1/2m_1} |E_1 - 2E_2 + E_3| \, dt \leq 5Mm^{-\alpha_1}$. Plugging this bound into (13.1), we have

$$(13.2) \qquad S_\ell \leq \left[\frac{M^2(n-m_0)}{m_1}\right] m_1^{-2\alpha_1} \qquad \text{for } \ell = 2, \ldots, m_1 - 1.$$

By definition, $r_0 = r_{m_1} = 0$ and the only bound available for the first and last term is $|r_\ell| \le M m_1^{-1}$, thus

$$(13.3) \qquad \max(S_1, S_{m_1}) \le \left[\frac{M^2(n - m_0)}{m_1}\right] m_1^{-2}.$$

DEPARTMENT OF STATISTICS AND APPLIED PROBABILITY
UNIVERSITY OF CALIFORNIA
SANTA BARBARA, CALIFORNIA 93106-3110
USA
E-MAIL: carter@pstat.ucsb.edu